\documentclass[11pt]{article}

\usepackage[]{graphicx,float,latexsym,times}
\usepackage{amsfonts,amstext,amsmath,amssymb,amsthm,dsfont}
 
\usepackage{enumerate}
\usepackage{setspace}
 
\RequirePackage[OT1]{fontenc}
\RequirePackage{xcolor}
\RequirePackage[numbers]{natbib}
\RequirePackage[colorlinks,citecolor=blue,urlcolor=blue]{hyperref}

\usepackage{setspace}

\usepackage{epsfig}
\usepackage{graphicx}
\usepackage{wrapfig}
\usepackage{amssymb}
\usepackage{amsfonts}
\usepackage{amsmath}
\usepackage{amsthm}
\usepackage{enumerate}
\usepackage{multicol}
\usepackage{bbm} 

\newtheorem{theorem}{Theorem}[section]
\newtheorem{proposition}[theorem]{Proposition}
\newtheorem{lemma}[theorem]{Lemma}

\newtheorem{remark}{Remark}[section]
\newtheorem{corollary}[theorem]{Corollary}
\newcommand\fdem{\hfill $\Box$}
\newcommand\cA{{\cal A}}
\newcommand\cC{{\cal C}}

\newcommand\cG{{\cal G}}

\newcommand\cD{{\cal D}}
\newcommand\cQ{{\cal Q}}

\newcommand\cR{{\cal R}}

\newcommand\ve{\varepsilon}
\newcommand\ov{\overline}

\def\bbr{{\mathbb R}}

\def\text#1{\hbox{#1}}
\def\proof{{\noindent \bf Proof. }}

\def\E{{\bf E}}

\def\P{{\bf P}}
\def\C{{\bf C}}

\def\D{{\bf D}}
\def\H{{\bf H}}

\def\U{{\bf U}}
\def\g{{\bf g}}
\def\l{{\bf l}}
\def\L{{\bf L}}

\newcommand{\wh}{\widehat}
\newcommand{\wt}{\widetilde}

\newcommand\Er{\mbox{Err}}

\newcommand\Tr{\mbox{Tr}}

\def\Chi{{\bf 1}}
\def\d{\mathrm{d}}
\def\build #1_#2{\mathrel{\mathop{\kern 0pt #1}\limits_\zs{#2}}}
\newcommand{\zs}[1]{{\mathchoice{#1}{#1}{\lower.25ex\hbox{$\scriptstyle#1$}}
{\lower0.25ex\hbox{$\scriptscriptstyle#1$}}}}
\numberwithin{equation}{section}

\begin{document}
\title{
Robust adaptive efficient  estimation for a semi-Markov  continuous
time regression from discrete data
\thanks{
This research was supported by RSF, project no 20-61-47043
 (National Research Tomsk State University , Russia).
}
}

\author{Vlad Stefan Barbu\thanks{
Laboratoire de Math\'ematiques Rapha\"el Salem,
 UMR 6085 CNRS-Universit\'e de Rouen Normandie,  France, e-mail: barbu@univ-rouen.fr},
Slim Beltaief\thanks{
Laboratoire de Math\'ematiques Rapha\"el Salem,
 UMR 6085 CNRS-Universit\'e de Rouen Normandie,  France, e-mail: slim.beltaief1@univ-rouen.fr
 }
 and 
 Serguei Pergamenshchikov\thanks{Laboratoire de Math\'ematiques Rapha\"el Salem,
 UMR 6085 CNRS-Universit\'e de Rouen Normandie,  France
and
International Laboratory of Statistics of Stochastic Processes and
Quantitative Finance of National Research Tomsk State University, Russia,
 e-mail:
Serge.Pergamenshchikov@univ-rouen.fr}}

\date{}
\maketitle

\begin{abstract}
In this article we consider the nonparametric robust estimation problem for regression models in continuous time with semi-Markov noises observed in discrete time moments.
 An adaptive model selection procedure   is proposed. A sharp 
non-asymptotic oracle inequality for the robust risks is obtained. We obtain sufficient conditions
 on the frequency observations under which the robust efficiency is shown.
 It turns out that for the semi-Markov models  the
 robust
  minimax convergence  rate may be faster or slower  than the classical one.
 \end{abstract}

\vspace*{5mm}

\vspace*{5mm}

\newpage

\section{Introduction}\label{sec:In}

In this paper we consider the semi-Markov regression model in continuous time introduced in \cite{VladSlimSerge2016}, i.e.
\begin{equation}\label{sec:In.1}
 \d\,y_\zs{t} = S(t)\d\,t + \d\,\xi_\zs{t}\,,\quad
 0\le t \le n\,,
\end{equation}
where $S(\cdot)$ is an unknown $1$-periodic function defined on $\bbr$ with values on $\bbr$, $(\xi_\zs{t})_\zs{t\ge 0}$
is the unobserved noise  process defined through a certain semi-Markov process in Section \ref{sec:Cns}.

Our problem in the present paper is to estimate the unknown function $S$ in the model \eqref{sec:In.1} on the basis of observations

\begin{equation}\label{sec:In.02}
(y_\zs{t_j})_\zs{0\le j\le np},\;\;\ t_j=j\Delta,\;\;\; \Delta=\frac{1}{p},
\end{equation}
where the integer $p\ge 1$ is the observation frequency. Firstly, this problem was considered in the framework ``signal+white noise'' (see, for example, \cite{IbragimovKhasminskii1981} or \cite{Pinsker1981}). 
Later, to introduce a dependence in the continuous time regression model
in \cite{KonevPergamenshchikov2003},
\cite{HopfnerKutoyants2009}, \cite{HopfnerKutoyants2010}
\cite{KonevPergamenshchikov2010}, 
 the Ornstein-Uhlenbeck processes has been  used
to model the ``color noise''. Moreover, in order to introduce the dependence and
the jumps in the regression model \eqref{sec:In.1},  the papers
\cite{KonevPergamenshchikov2012} and \cite{KonevPergamenshchikov2015}
use the non Gaussian Ornstein-Uhlenbeck processes defined in
\cite{BarndorffNielsenShephard2001}. The problem in all these papers is that the introduced
Ornstein-Uhlenbeck type of dependence decreases  with a geometric rate. So, asymptotically
when the duration of observations goes to infinity,
we obtain the same
``signal+white noise'' model
very quick.
To keep the dependence for sufficiently large duration of observations, 
in \cite{VladSlimSerge2016} it was proposed the model \eqref{sec:In.1} with a semi-Markov component in the 
jumps of the noise process $(\xi_\zs{t})_\zs{t\ge 0}$.

 The main goal of this paper is to develop adaptive robust method from \cite{VladSlimSerge2016}, that was based on continuous observations, to the 
 estimation problem based on discrete observations given in \eqref{sec:In.02}.
In this paper we use  quadratic risk defined as
\begin{equation}\label{sec:In.4}
\cR_\zs{Q}(\wt{S}_\zs{n},S)=
\E_\zs{Q,S}\,\|\wt{S}_\zs{n}-S\|^{2}\,,
\end{equation}
where $\wt{S}_\zs{n}(\cdot)$ is some estimate (i.e. any periodical function measurable with respect to the observations 
$\sigma\{y_\zs{t_\zs{0}},\ldots y_\zs{t_\zs{n}}\}$),
 $\|f\|^{2}=\int^{1}_\zs{0}\,f^{2}(s)\d s$ and
$\E_\zs{Q,S}$ is the expectation with respect to the distribution $\P_\zs{Q,S}$
 of the process  \eqref{sec:In.1} corresponding to the unknown noise distribution $Q$ in 
  the Skorokhod space $\cD[0,n]$.
We assume that this distribution belongs to some distribution family $\cQ_\zs{n}$ specified  in Section \ref{sec:Cns}.
 
To study the properties of the estimators uniformly over the noise distribution
(what is really needed in practice),  we use the robust risk defined as
\begin{equation}\label{sec:In.6}
\cR^{*}_\zs{n}(\wt{S}_\zs{n},S)=\sup_\zs{Q\in\cQ_\zs{n}}\,
\cR_\zs{Q}(\wt{S}_\zs{n},S)\,.
\end{equation}

 Thus the goal of this paper is to develop a robust efficient model selection  method
 based on the observations
 \eqref{sec:In.02}
for the model \eqref{sec:In.1} with the semi-Markov components in the jumps of the  noise $(\xi_\zs{t})_\zs{t\ge 0}$.
 We use the approach proposed by Konev and Pergamenshchikov in
  \cite{KonevPergamenshchikov2015}
 for continuous-time regression models observed in the discrete time moments.
   Unfortunately, we cannot use directly this method
  for semi-Markov regression models, since their tool essentially uses
  the fact that the Ornstein-Uhlenbeck dependence decreases with  geometrical rate
  and obtain sufficiently quickly the ``white noise'' case.  In the present paper, in order to obtain the sharp non-asymptotic oracle inequalities, we use the
    renewal  methods from \cite{VladSlimSerge2016} developed for the model \eqref{sec:In.1}.
 As  a consequence, we can obtain
  the constructive sufficient conditions that provide the
   robust efficiency for proposed model selection procedures.

 The rest
of the paper is organized as follows. In Section~\ref{sec:Cns}
we state the main conditions under which we consider the model  \eqref{sec:In.1}.
In Section~\ref{sec:Mo}
 we
construct the model selection procedure on the basis of
weighted least squares estimates, here we also specify the set
 of admissible weight sequences in the model selection procedure.
In Section~\ref{sec:Mrs}
we state the main results in the form of oracle inequalities for the quadratic
 risk  and the robust risk. In Section \ref{sec:Prsm} we study some properties of the regression model \eqref{sec:In.1}. Section~\ref{sec:Siml} is devoted to some numerical results.
 In section~\ref{sec:SI} we study some properties of the stochastic integral. Section~\ref{sec:Pr}  gives the proofs of the oracle inequalities
for the regression model \eqref{sec:In.1} with the noises
introduced in Section~\ref{sec:Cns}. Some auxiliary are given in an Appendix.

\bigskip

\section{Main conditions}\label{sec:Cns}

First, we assume that 
the noise process $(\xi_\zs{t})_\zs{t\ge\, 0}$ in the model \eqref{sec:In.1} is defined  as
\begin{equation}\label{sec:Ex.1}
 \xi_\zs{t} =\varrho_{1} w_\zs{t}+  \varrho_{2} L_\zs{t} + \varrho_{3} z_\zs{t}\,,
\end{equation}
where $\varrho_\zs{1}$, $\varrho_\zs{2}$ and $\varrho_\zs{3}$ are unknown coefficients,
$(w_\zs{t})_\zs{t\ge\,0}$ is a standard Brownian motion, $(L_\zs{t})_\zs{t\ge\,0}$ is a jump L\'evy process defined as
 \begin{equation}\label{sec:Mcs.01}         
L_\zs{t}=
\int^{t}_\zs{0}\int_\zs{\bbr_\zs{*}}x (\mu(\d s,\d x)
-\wt{\mu}(\d s,\d x))
\,,
\end{equation}
 $\mu(\d s\,\d x)$ is the jump measure with deterministic
compensator $\wt{\mu}(\d s\,\d x)=\d s\Pi(\d x)$, 
$\Pi(\cdot)$  is the Levy measure on $\bbr_\zs{*}=\bbr\setminus\{0\}$
(see, for example
\cite{JacodShiryaev2002, ContTankov2004} for details), with
\begin{equation}\label{sec:Mcs.2}
\Pi(x^{2})=1
\quad\mbox{and}\quad
\Pi(x^{8})
\,<\,\infty\,.
\end{equation}
Here we use the usual notations for $\Pi(\vert x\vert^{m})=\int_\zs{\bbr}\,\vert z\vert^{m}\,\Pi(\d z)$.
Note that  $\Pi(\vert x\vert)$ may be equal to $+\infty$.
In this paper we assume that the ``dependent part'' in the noise 
 \eqref{sec:Ex.1} is modeled  by the semi- Markov process $(z_\zs{t})_\zs{t\ge\, 0} $ defined as
 \begin{equation}\label{sec:Ex.2}
 z_\zs{t} = \sum_\zs{i=1}^{N_\zs{t}} Y_\zs{i},
\end{equation}
where $(Y_\zs{i})_\zs{i\ge\, 1}$ is a sequence of random variables verifies the following Condition\\ 
%
$\C_\zs{*})$
$$ \E_\zs{Q}\,Y_\zs{i}^2=1,  \E_\zs{Q}\,Y_\zs{i_1}\, Y_\zs{i_2} =0 \quad\mbox{and}\quad \sup_\zs{i \geq 1}\E_\zs{Q}\,Y^4_\zs{i}<\infty
\,.$$
$\quad\mbox{and}\quad$
$$\forall\, i_1< i_2< i_3< i_4,\,  p_i \in\{0,1,2,3,4\} \quad\mbox{such that}\quad  \sum_{i=1}^4 p_i =4 $$
$$ \E_\zs{Q}\,Y^{p_1}_\zs{i_1}\, Y^{p_2}_\zs{i_2}\, Y^{p_3}_\zs{i_3}\, Y^{p_4}_\zs{i_4}=0\,\,If\,\, \exists\,\,  p_i - odd  $$
Now let us give some examples:\\\\
\bigskip
$1)$  $(Y_\zs{j})_\zs{j\geq\,1}$  independent  with 
$$ \E_\zs{Q}\,Y_\zs{j}=0, \E_\zs{Q}\,Y_\zs{j}^2=1 \quad\mbox{and}\quad \sup_\zs{j \geq 1}\E_\zs{Q}\,Y^4_\zs{j}<\infty.$$
$2)$  $\exists\,d>1$ such that  the random vector $(Y_\zs{1},...,Y_\zs{d})$ has  a  spherically symmetric distribution
with $\E_\zs{Q}\,Y_\zs{1}^2=1,\, \E_\zs{Q}\,Y_\zs{1}^4 <  \infty $
$\quad\mbox{and}\quad$  $(Y_\zs{j})_\zs{j\geq\,d}$ are such that 
$$ \E_\zs{Q}\,Y_\zs{j}^2=1  \quad\mbox{and}\quad  \sup_\zs{j\geq d}\E_\zs{Q}\,Y^4_\zs{j}<\infty.$$
$3)$ $(Y_\zs{j})_\zs{j\geq\,1}$  such that $\forall n\geq 1$,  $(Y_1,...,Y_n)$  is a Gaussian Mixture.\\\\

\bigskip
Here $N_\zs{t}$ is a general counting process (see, for example, \cite{Mikosch2004})
defined as
\begin{equation}\label{sec:Ex.4}
N_\zs{t} = \sum_\zs{k=1}^{\infty} 1_\zs{\{T_\zs{k} \le t\}}
\quad\mbox{and}\quad
T_\zs{k}=\sum_\zs{l=1}^k\, \tau_\zs{l}\,,
\end{equation}
with $(\tau_\zs{l})_\zs{l\ge\,1}$ an i.i.d. sequence of positive integrated
 random variables with the distribution $\eta$ and mean $\check{\tau}=\E\,\tau_\zs{1}>0$.
 We assume that the processes 
$(N_\zs{t})_\zs{t\ge 0}$ and  $(Y_\zs{i})_{i\ge\, 1}$ are independent between them and are also independent of $(L_\zs{t})_\zs{t\ge 0}$.
Note that the process $(z_\zs{t})_\zs{t\ge\, 0}$ is a special case of a semi-Markov process (see, e.g., \cite{BarbuLimnios2008} and \cite{LO}).
\begin{remark} \label{Re.sec.Ex.0}
It should be noted that, if $\tau_\zs{j}$ is an Exponential random variable, i.e. $g$ is the Exponential density, then $(N_\zs{t})_\zs{t\ge 0}$
is a Poisson process and, in this case, $(\xi_\zs{t})_\zs{t\ge 0}$ is a L\'evy process for which this model is studied
in   \cite{KonevPergamenshchikov2009a}, \cite{KonevPergamenshchikov2009b} and \cite{KonevPergamenshchikov2012}.
But, in the general case when the process \eqref{sec:Ex.2} is not a L\'evy process, this process has a memory and cannot be treated
in the framework of semi-martingales with independent increments. One needs to 
develop a new tool based on the renewal theory arguments. 
\end{remark}
\begin{remark} \label{Re.sec.Ex.0}
Note that the noise in our model is the sum of the L\'evy process given by  $(\varrho_{1} w_\zs{t}+\varrho_{2} L_\zs{t})_\zs{t\ge\,0}$  and the semi-Markov process $(z_\zs{t})_\zs{t\ge\, 0} $  in order to include dependent observations.
\end{remark}
Let us denote by $\rho$  the density of
 the renewal measure 
$\check{\eta}$ defined as
\begin{equation}\label{sec:Cns.1}
\check{\eta}
=\sum^{\infty}_\zs{l=1}\,\eta^{(l)}
\,,
\end{equation}
where $\eta^{(l)}$ is the $l$th convolution power of the measure $\eta$. 
As to the parameters in  \eqref{sec:Ex.1},
we assume that 
\begin{equation}\label{sec:Ex.5}
\varsigma_\zs{*}\le 
\sigma_\zs{Q}
\le \varsigma^{*}\,,
\end{equation}
where $\sigma_\zs{Q}=\varrho_{1}^{2}+\varrho_{2}^{2}+ \varrho_{3}^{2}/\check{\tau}$, 
the    unknown 
bounds $0<\varsigma_\zs{*}\le \varsigma^{*}$ are functions of $n$, i.e.  $\varsigma_\zs{*}=\varsigma_\zs{*}(n)$ and  $\varsigma^{*}=\varsigma^{*}(n)$, such that
for any $\check{\epsilon}>0$
\begin{equation}\label{sec:Mrs.5-1}
\lim_\zs{n\to\infty}n^{\check{\epsilon}}\,\varsigma_\zs{*}(n)=+\infty
\quad\mbox{and}\quad
\lim_\zs{n\to\infty}\,\frac{\varsigma^{*}(n)}{n^{\check{\epsilon}}}=0\,.
\end{equation}
We denote by $\cQ_\zs{n}$ the family of all distributions of the process \eqref{sec:Ex.1}
in $\D[0,n]$ satisfying the properties \eqref{sec:Ex.5} -- \eqref{sec:Mrs.5-1}.

\begin{remark} \label{Re.sec.Ex.1__00}
As we will see later,  the parameter $\sigma_\zs{Q}$ is the limit of the Fourier transform of the noise process \eqref{sec:Ex.1}.
 Such a limit is called variance proxy
 (see \cite{KonevPergamenshchikov2012}).
\end{remark}

We assume that the distribution $\eta$  has a density $g$ that satisfies the following conditions.

\bigskip

$\H_\zs{1}$) {\em Assume that, for any $x\in\bbr,$ 
there exist the finite limits
$$
g(x-)=\lim_\zs{z\to x-}g(z)
\quad\mbox{and}\quad
g(x+)=\lim_\zs{z\to x+}g(z)
$$
and, for any $K>0,$
there exists $\delta=\delta(K)>0$
for which
\begin{equation}\label{sec:A.8-00}
\sup_\zs{\vert x\vert\le K}\,
\int^{\delta}_\zs{0}\,
\frac{
\vert 
g(x+t)+g(x-t)-g(x+)-g(x-)
\vert
}{t}
\d t
\,<\,\infty.
\end{equation}
}

\bigskip

$\H_\zs{2}$) {\em  For any $\gamma>0,$
$$
\sup_\zs{z\ge 0}\,z^{\gamma}\vert 2g(z) -g(z-)-g(z+) \vert\,<\,\infty\,.
$$
}

$\H_\zs{3}$) {\em  There exists $\beta>0$ such that $\int_\zs{\bbr}\,e^{\beta x}\,g(x)\,\d x<\infty.$
}

\bigskip
$\mathbf{Example}$

\bigskip
Let $(N_\zs{t})_\zs{t\ge 0}$ be a truncated Fractional Poisson process defined as in  \eqref{sec:Ex.4} where $(\tau_\zs{k})_\zs{k\ge\,1}$ is an i.i.d. sequence with the 
Mittag-Leffler distribution (see, for example, \cite{BiardSaussereaur2014}  given by 
$$
\mathbf{P}(\tau_\zs{1}> t) = \rho_\zs{H}(- \lambda t^h) 
$$
for $\lambda>0$ and $0< H < 1$, where
$$
\rho_\zs{H}(z) = \sum_{k=0}^{\infty}  \frac{z^k}{\Gamma(1+H k)}. 
$$
$\tau_\zs{1}$ has the same distribution as $X_0 \wedge X_1$, $X_0$ is a Mittag-Leffler waiting time with distribution $ \rho_\zs{H}$ and $X_1$ is an exponential random variable with intensity $  \check{\lambda} = 1-H $ and density $  \check{\rho} $ ,  in this case 
$\tau_\zs{1}$ is light lailed , i.e
$$ 
\exists r>0,\;\;
\E\left( e^{r \tau_\zs{1}}\right)  < +\infty,
$$
it's easy to check that  the density $g$  of the variables  $(\tau_\zs{k})_\zs{k\ge\,1}$  is given by 
$$
g(z) = \rho_\zs{H}(z) P(X_1>z) + \check{\rho} (z) P(X_0>z)
$$
and  verifies conditions  $\H_\zs{1})$--$\H_\zs{3})$

\bigskip
\begin{remark} \label{Re.sec.Ex.201}
In the case where $\tau_\zs{1} = X_0$ conditions  $\H_\zs{1})$--$\H_\zs{3})$ don't hold. Indeed, the variable $ \tau_\zs{1}$ has heavy tails and infinite mean, i.e
$$ 
\forall r>0,\; \; \ \E\left( e^{r \tau_\zs{1 }} \right) = +\infty.
$$
\end{remark}

%
%

\bigskip
\begin{remark} \label{Re.sec.Ex.2}
It should be noted that Condition $\H_\zs{3})$ means that there exists an exponential moment for the random variable $(\tau_\zs{j})_\zs{j\ge 1}$, i.e. these random variables are not too large.  This is a natural 
constraint since these random variables define the intervals between jumps, i.e. the jump frequency. 
 So, to study the influence of the jumps in the model \eqref{sec:In.1} one needs to consider the noise process \eqref{sec:Ex.1} with ``small'' intervals between jumps or
 large jump frequency. 
\end{remark}

For the next condition we need the Fourier transform for any function $f:\bbr\to\bbr$ from $\L_\zs{1}(\bbr)$ defined 
by
\begin{equation}\label{sec:Rtl.06-0}
\wh{f}(\theta)=\frac{1}{2\pi}\,\int_\zs{\bbr}\,e^{i\theta x}\,f(x)\,\d x\,.
\end{equation}

\bigskip

$\H_\zs{4}$) {\em  There exists $t^{*}>0$ such that  the function $\wh{g}(\theta-it)$ 
belongs to $\L_\zs{1}(\bbr)$ for any $0\le t\le t^{*}$.
}

It is clear that Conditions $\H_\zs{1})$--$\H_\zs{4})$ hold true for
any continuously differentiable function $g$ having an exponential moment, for example, for
 the $\Gamma$ density.

It should be noted that in view of Proposition 5.2 from \cite{VladSlimSerge2016},
 Conditions $\H_\zs{1})$--$\H_\zs{4})$ imply 
\begin{equation}\label{sec:Mrs.1-0-1}
\Vert\Upsilon\Vert_\zs{1}=\int^{+\infty}_\zs{0}\,\vert\Upsilon(x)\vert\,\d x
\,<\infty\,,
\end{equation}
where $\Upsilon(x)=\rho(x)-1/\check{\tau}$.

\section{Model selection}\label{sec:Mo}

In this section we construct a model selection procedure for estimating the unknown function $S$ given in\eqref{sec:In.1} starting from the discrete-time observations \eqref{sec:In.02} and we establish the oracle inequality for the associated risk. To this end, note that
for any function $f:[0,n] \to\bbr$ from $\L_\zs{2}[0,n]$, the integral
\begin{equation}\label{sec:In.2}
I_\zs{n}(f)=\int_\zs{0}^{n} f(s) \d\xi_\zs{s}
\end{equation}
is well defined, with $\E_\zs{Q}\,I_n(f)=0$. Moreover, as it is shown in Lemma \ref{Le.sec:Smp.1} 
under the conditions $\H_\zs{1})$--$\H_\zs{4}),$
\begin{equation}\label{sec:In.3}
\E_\zs{Q}\,I^{2}_n(f) \le \varkappa_\zs{Q}\,\int_\zs{0}^{n} f^2_\zs{s} \d\,s
\quad\mbox{and}\quad
\varkappa_\zs{Q}=\bar{\varrho}+\varrho_{3}^2\,\vert\rho\vert_\zs{*}
\end{equation}
where $\bar{\varrho}=\varrho_{1}^{2}+\varrho_{2}^{2}$ and $\vert\rho\vert_\zs{*}=\sup_\zs{t\ge 0}\vert\rho(t)\vert<\infty$.
\noindent 

\noindent 
In this paper we will use the trigonometric basis  $(\phi_\zs{j})_\zs{j\ge\, 1}$ in $\L_\zs{2}[0,1]$ defined as
\begin{equation}\label{sec:In.5}
\phi_\zs{1} = 1\,,\quad \phi_\zs{j}(x)= \sqrt 2 \Tr_\zs{j}(2\pi[j/2]x)\,,\quad j\ge\,2\,,
\end{equation}
where the function $\Tr_\zs{j}(x)= \cos(x)$ for even $j$ and $\Tr_\zs{j}(x)= \sin(x)$ for odd $j$, $[x]$ denotes the integer part of $x$. By making use of this basis we consider the discrete Fourier transformation of $S$
\begin{equation}\label{sec:In.06_Fou}
S(t)=  \sum_\zs{j=1}^{p}\,\theta_\zs{j,p}\, \phi_\zs{j}(t),\;\;\; t\in{\{t_1,...,t_p}\},
\end{equation}
where  the Fourier coefficients are defined by
\begin{equation}\label{sec:In.6_Fou}
\theta_\zs{j,p}=(S,\phi_\zs{j})_\zs{p} = \frac{1}{p} \sum_\zs{i=1}^{p} S(t_i) \phi_\zs{j}(t_i).
\end{equation}
In the sequel the corresponding norm will be denoted by $\Vert x\Vert^{2}_\zs{p}=(x,x)_\zs{p}$.
These Fourier coefficients $\theta_\zs{j,p}$ can be estimated by
\begin{equation}\label{sec:In.7}
\wh{\theta}_\zs{j,p}= \frac{1}{n} \int_\zs{0}^{n}  \Psi_\zs{j,p}(t) \d\,y_\zs{t},
\quad\mbox{and}\quad
\Psi_\zs{j,p}(t) = \sum_\zs{l=1}^{np} \phi_\zs{j}(t_l)1_\zs{\{t_{l-1} < t \leq t_l \}}
\,.
\end{equation}

Let us note that the system of the functions $(\Psi_\zs{j,p})_{1 \leq j \leq p }$ is orthonormal in $\L_\zs{2}[0,1]$ because
$$
 \int_\zs{0}^{1} \Psi_\zs{j,p}(t) \Psi_\zs{i,p}(t) \d\,t\ = 
 (\phi_\zs{j},\phi_\zs{i})_\zs{p}
 =1_\zs{\{i=j \}}
\,.
$$
In the sequel we need the Fourier coefficients of the function $S$ with respect to the new basis $(\Psi_\zs{j,p})_{1 \leq j \leq p }.$
 These coefficients can be written as 
\begin{equation}\label{sec:In.08}
\overline{\theta}_\zs{j,p} = \int_\zs{0}^{1} S(t) \Psi_\zs{i,p}(t) \d\,t\ = \theta_\zs{j,p} + h_{j,p},
\end{equation}
where
$$
h_{j,p}(S)= \sum_\zs{l=1}^{p} \int_\zs{t_{l-1}}^{t_l} \phi_\zs{j}(t_l) (S(t)-S(t_l)) \d\,t
\,.
$$

From \eqref{sec:In.1} it follows directly that these Fourier coefficients satisfy the equation
\begin{equation}\label{sec:In.8}
\wh{\theta}_\zs{j,p}= \overline{\theta}_\zs{j,p} + \frac{1}{\sqrt n}\xi_\zs{j,p},
\quad\mbox{where}\quad
\xi_\zs{j,p}= \frac{1}{\sqrt n} I_\zs{n}(\Psi_\zs{j,p})
\,.
\end{equation}


\noindent
For any $0\le t\le 1$
we estimate the function $S$ by the weighted least squares estimator
\begin{equation}\label{sec:Mo.1}
\wh{S}_\zs{\lambda} (t) = \sum_\zs{j=1}^{n} \lambda(j) \wh{\theta}_\zs{j,p} \Psi_\zs{j,p}(t)\,,
\end{equation}
where the weight vector $\lambda=(\lambda(1),.....,\lambda(n))$ belongs to some finite set $\Lambda$ from $[0,1]^n$, $\wh{\theta}_\zs{j,n}$ was defined in \eqref{sec:In.7} and $\phi_\zs{j}$ in \eqref{sec:In.5}.  Now let us consider
 \begin{equation}\label{sec:Mo.2}
\nu=\#(\Lambda)
\quad\mbox{and}\quad
\vert\Lambda\vert_\zs{*}= \max_\zs{\lambda\in\Lambda}\,L(\lambda)
\,,
\end{equation}
where $\#(\Lambda)$ is the cardinal number of  $\Lambda$ and $L(\lambda)=\sum^{n}_\zs{j=1}\lambda(j)$. In the sequel we assume that
$\vert\Lambda\vert_\zs{*}\ge 1$ and $\lambda(j)=0$ for $j\ge p$.

In order to find a proper weight sequence $\lambda$ in the set $\Lambda,$ one needs to specify a cost function. When choosing an appropriate cost function, one can use the following argument. Let as consider the empirical squared error 
\begin{equation}\label{sec:def-err.1}
\Er(\lambda) = \Vert \wh{S}_\lambda-S\Vert^2\,,
\end{equation}
which in our case is equal to
\begin{equation}\label{sec:Mo.3}
\Er(\lambda) = \sum_\zs{j=1}^{n} \lambda^2(j) \wh{\theta}^2_\zs{j,p} -2 \sum_\zs{j=1}^{n} \lambda(j) \wh{\theta}_\zs{j,p}\overline{\theta}_\zs{j,p}+
\Vert S\Vert^{2}
\,.
\end{equation}
Since the Fourier coefficients $(\theta_\zs{j})_\zs{j\ge\,1}$ are unknown, the weight coefficients $(\lambda(j))_\zs{1 \leq j\leq p}$ cannot be determined by minimizing this quality. To circumvent this difficulty, one needs to replace the terms $\wh{\theta}_\zs{j,p}\overline{\theta}_\zs{j,p}$ by
their estimators $\wt{\theta}_\zs{j,p}$. Let us set
\begin{equation}\label{sec:Mo.4}
\wt{\theta}_\zs{j,p} = \wh{\theta}^2_\zs{j,p} - \frac{  \wh{\sigma}_\zs{n}}{n}.
\end{equation}
Here $\wh{\sigma}_\zs{n}$ is an estimate for the proxy variance $\sigma_\zs{Q}$
defined in
\eqref{sec:Ex.5}. For, example, we can take it as
\begin{equation}\label{sec:Mo.4-1-31-3}
\wh{\sigma}_\zs{n}=\frac{n}{\check{p}}
\sum^{\check{p}}_\zs{j=l}\,\wh{\theta}^2_\zs{j,p}
\quad\mbox{and}\quad
\check{p}=\min(p,n),
\end{equation}
where $ l=[\sqrt{n}]$ and we set $\wh{\sigma}_\zs{n} = 0 \; for \; l > p.$
For this change in the empirical squared error, one has to pay some penalty. Thus we obtain the cost function of the form
\begin{equation}\label{sec:Mo.5}
J_n(\lambda)=\sum_\zs{j=1}^{n} \lambda^2(j) \wh{\theta}^2_\zs{j,n} -2 \sum_\zs{j=1}^{n} \lambda(j)\wt{\theta}_\zs{j,n} + \delta\,P_\zs{n}(\lambda),
\end{equation}
where $\delta>0$ is some  threshold which will be specified later and the penalty term is
\begin{equation}\label{sec:Mo.6}
P_n(\lambda)= \frac{  \wh{\sigma}_\zs{n} |\lambda|^2}{n}.
\end{equation}
Minimizing the cost function, that is
\begin{equation}\label{sec:Mo.8}
\hat\lambda= \mbox{argmin}_\zs{\lambda\in\Lambda} J_n(\lambda),
\end{equation}
and substituting the obtained weight coefficients $\hat\lambda$ in \eqref{sec:Mo.1}, lead to the model selection procedure
\begin{equation}\label{sec:Mo.9}
\wh{S}_* = \wh{S}_\zs{\hat \lambda}.
\end{equation}
We recall that the set $\Lambda$ is finite, so $\hat \lambda$ exists. In the case when $\hat \lambda$ is not unique we take one of them.

\section{Main results}\label{sec:Mrs}

\subsection{Oracle inequalities}

First we define the following constant which will be used 
to describe the rest term 
in the oracle inequalities. We set

\begin{equation}\label{sec:Mrs.1-0}
\g_\zs{n,p}
=1+\vert\Lambda\vert_\zs{*}\left(\frac{\sqrt{n}}{\check{p}}+ \frac{1}{\sqrt{\check{p}}}\right)
\,.
\end{equation}

\noindent
Firstly,  we obtain the non asymptotic oracle inequality for the model selection procedure \eqref{sec:Mo.9}.

\begin{theorem}\label{Th.sec:Mrs.0}
Assume that 
Conditions $\H_\zs{1})$--$\H_\zs{4})$
hold true. Then, 
there exists some constant $\l^{*}>0$ such that, for any 
noise distribution
$Q$,  the weight vector set $\Lambda$, for any periodic function $S$
for any $n\ge\, 1 $, $p\ge 3$ and $ 0 <\delta \leq 1/6$, 
 the procedure \eqref{sec:Mo.9} 
satisfies the following oracle inequality
\begin{align}\nonumber
\mathcal{R}_\zs{Q}(\wh{S}_*,S)\leq&\frac{1+3\delta}{1-3\delta} \min_\zs{\lambda\in\Lambda} 
\mathcal{R}_\zs{Q}(\wh{S}_\lambda,S)
\\[2mm]  \label{sec:Mrs.1--00}
&
+
\l^{*}\frac{\nu}{\delta n}
\left(
\sigma_\zs{Q}
+\vert\Lambda\vert_\zs{*}\,
\E_\zs{Q}\vert \wh{\sigma}_\zs{n}
-
\sigma_\zs{Q}\vert
\right)
\,.
\end{align}
\end{theorem}

  \begin{corollary}\label{Co.sec:Mrs.0}
Assume that 
Conditions $\H_\zs{1})$--$\H_\zs{4})$
hold true and that the proxy variance $\sigma_\zs{Q}$ is known. 
Then there exists some constant $\l^{*}>0$ such that for any 
noise distribution
$Q$,  the weight vectors set $\Lambda$, for any periodic function $S$
for any $n\ge\, 1 $, $p\ge 3$ and $ 0 <\delta \leq 1/6$, 
 the procedure \eqref{sec:Mo.9} with $\wh{\sigma}_\zs{n}=\sigma_\zs{Q}$, 
satisfies the following oracle inequality
\begin{equation}\label{sec:Mrs.1_cor_1}
\mathcal{R}_\zs{Q}(\wh{S}_*,S)\leq\frac{1+3\delta}{1-3\delta} \min_\zs{\lambda\in\Lambda} 
\mathcal{R}_\zs{Q}(\wh{S}_\lambda,S)
+
\l^{*}\frac{\sigma_\zs{Q}\nu}{\delta n}
\,.
\end{equation}
\end{corollary}

\noindent 
Now we study the model selection 
procedure \eqref{sec:Mo.9} using the proxy estimate
  \eqref{sec:Mo.4-1-31-3}.
\begin{theorem}\label{Th.sec:Mrs.1-10}
Assume that the function $S$ is continuously differentiable 
 and that 
Conditions $\H_\zs{1})$--$\H_\zs{4})$
hold true. 
Then 
there exists some constant $\l^{*}>0$ such that for any 
noise distribution
$Q$,  the weight vectors set $\Lambda$, for any periodic function $S$
for any $n\ge\, 1 $, $p\ge 3$ and $ 0 <\delta \leq 1/6$, 
 the procedure \eqref{sec:Mo.9} 
satisfies the following oracle inequality
\begin{align}\nonumber
\mathcal{R}_\zs{Q}(\wh{S}_*,S)\leq&\frac{1+3\delta}{1-3\delta} \min_\zs{\lambda\in\Lambda} 
\mathcal{R}_\zs{Q}(\wh{S}_\lambda,S)\\[2mm]  \label{sec:Mrs.1}
&+
\l^{*}\frac{\nu}{\delta n}
(1+\sigma_\zs{Q})^{3}
\,
\left(
1+\Vert \dot{S}\Vert^{2}
\right)\g_\zs{n,p}
\,.
\end{align}
\end{theorem}

\noindent 
Let us study the robust risks
\eqref{sec:In.6} for the procedure \eqref{sec:Mo.9}.  
In this case this family  consists  of all distributions on the Skorokhod space $\cD[0,n]$ of the process
\eqref{sec:Ex.1}  with the parameters satisfying Conditions \eqref{sec:Ex.5}--\eqref{sec:Mrs.5-1}.

In order to obtain the efficiency property, we specify the weight coefficients 
$(\lambda(j))_\zs{1\le j\le n}$ in the procedure \eqref{sec:Mo.9}.
 Consider, for some fixed $0<\varepsilon<1,$ a numerical grid of the form
\begin{equation}\label{sec:Ga.0}
\cA=\{1,\ldots,k^*\}\times\{\varepsilon,\ldots,m\varepsilon\},
\end{equation}
where $m=[1/\ve^2]$. We assume that both
parameters $k^*\ge 1$ and $\varepsilon$ are functions of $n$, i.e.
$k^*=k^*(n)$ and $\ve=\ve(n)$, such that
\begin{equation}\label{sec:Ga.1}
\left\{
\begin{array}{ll}
&\lim_\zs{n\to\infty}\,k^*(n)=+\infty\,,
\quad\lim_\zs{n\to\infty}\,\dfrac{k^*(n)}{\ln n}=0\,,\\[6mm]
&
\lim_\zs{n\to\infty}\,\varepsilon(n)=0
\quad\mbox{and}\quad
\lim_\zs{n\to\infty}\,n^{\check{\delta}}\ve(n)\,=+\infty
\end{array}
\right.
\end{equation}
for any $\check{\delta}>0$. One can take, for example, for $n\ge 2$
\begin{equation}\label{sec:Ga.1-00}
\ve(n)=\frac{1}{ \ln n }
\quad\mbox{and}\quad
k^*(n)=k^{*}_\zs{0}+\sqrt{\ln n}\,,
\end{equation}
where $k^{*}_\zs{0}\ge 0$ is some fixed constant.
 For each $\alpha=(\beta, \l)\in\cA$, we introduce the weight
sequence
$$
\lambda_\zs{\alpha}=(\lambda_\zs{\alpha}(j))_\zs{1\le j\le p}
$$
with the elements
\begin{equation}\label{sec:Ga.2}
\lambda_\zs{\alpha}(j)=\Chi_\zs{\{1\le j<j_\zs{*}\}}+
\left(1-(j/\omega_\alpha)^\beta\right)\,
\Chi_\zs{\{ j_\zs{*}\le j\le \omega_\zs{\alpha}\}},
\end{equation}
where
$j_\zs{*}=1+\left[\ln\upsilon_\zs{n}\right]$, $\omega_\zs{\alpha}=(\d_\zs{\beta}\,\l\upsilon_\zs{n})^{1/(2\beta+1)}$,
$$
\d_\zs{\beta}=\frac{(\beta+1)(2\beta+1)}{\pi^{2\beta}\beta}
\quad\mbox{and}\quad
\upsilon_\zs{n}=n/\varsigma^{*}
\,.
$$
We remind that the threshold $\varsigma^{*}$ is introduced in the definition of the distribution family $\cQ_\zs{n}$ in
\eqref{sec:Ex.5}.
Now we define the set $\Lambda$ 
as
\begin{equation}\label{sec:Ga.3}
\Lambda\,=\,\{\lambda_\zs{\alpha}\,,\,\alpha\in\cA\}\,.
\end{equation}

These weight coefficients are
  used in \cite{KonevPergamenshchikov2012, KonevPergamenshchikov2015}
 for continuous time regression
models  to show the asymptotic efficiency.  Note also that in this case the cardinal of the set $\Lambda$ is  
\begin{equation}
\label{sec:Ga.1++1--1}
\nu=k^{*} m\,.
\end{equation}
Moreover,
taking into account that $\d_\zs{\beta}<1$ for $\beta\ge 1$ 
we obtain for the set \eqref{sec:Ga.3}
\begin{equation}
\label{sec:Ga.1++1--2}
 \vert \Lambda\vert_\zs{*}\,
 \le\,1+
\sup_\zs{\alpha\in\cA}
  \omega_\zs{\alpha}
\le 1+(\upsilon_\zs{n}/\ve )^{1/3}\,.
\end{equation}
Therefore, the last  condition in \eqref{sec:Ga.1} yields
$$
\lim_\zs{n\to\infty}\frac{\vert\Lambda\vert_\zs{*}}{n^{1/3+\check{\epsilon}}}=0
\quad\mbox{for any}\quad \check{\epsilon}>0
\,.
$$

Our goal is to bound asymptotically the term 
\eqref{sec:Mrs.1-0}
by any power of $n$. To this end, we assume the following condition
on the frequency of the observations.\\

$\H_\zs{5}$) {\em Assume that there exists $\check{\delta} > 0$ such that 
for any $n\ge 3$
\begin{equation}\label{sec:freq-cond}
p \ge n^{5/6}
\,.
\end{equation}
}

\noindent
Now, Theorem \ref{Th.sec:Mrs.1-10} implies the following oracle inequality.

\begin{theorem}\label{Th.sec:Mrs.2}
Assume that the unknown function $S$ is continuously differentiable. Moreover, assume that
Conditions $\H_\zs{1})$--$\H_\zs{5})$ 
hold true. Then, for the robust risks defined in  \eqref{sec:In.6} through the distribution family \eqref{sec:Ex.5}--\eqref{sec:Mrs.5-1}, the procedure \eqref{sec:Mo.9} with the coefficients
\eqref{sec:Ga.2}
for any $n\ge\, 1 $ and $ 0 <\delta <1/6$ 
satisfies the following oracle inequality  
\begin{equation}\label{sec:Mrs.6-25.3}
\cR^{*}(\wh{S}_*,S)\leq\frac{1+3\delta}{1-3\delta} \min_\zs{\lambda\in\Lambda}
 \cR^{*}(\wh{S}_\lambda,S)+
\frac{\U^{*}_\zs{n}(S)}{n\delta}
\,,
\end{equation}
where the sequence $\U^{*}_\zs{n}(S)>0$ is such that under condition \eqref{sec:Ga.1} 
for any $r>0$ and $\check{\delta}>0,$
\begin{equation}\label{sec:Mrs.7-25.3}
\lim_\zs{n\to\infty}\,
\sup_\zs{\|\dot{S}\| \le r}
\,
\frac{\U^{*}_\zs{n}(S)}{n^{\check{\delta}}}
=0
\,.
\end{equation}
\end{theorem}

\subsection{Robust asymptotic efficiency}

\noindent 
Now we study the asymptotically efficiency properties for the procedure
\eqref{sec:Mo.9}, \eqref{sec:Ga.2}
  with respect to the robust risks \eqref{sec:In.6} defined by the 
  distribution family \eqref{sec:Ex.5} -- \eqref{sec:Mrs.5-1}.  To this end, we assume that the unknown function $S$ in the model
 \eqref{sec:In.1} belongs to the Sobolev ball
\begin{equation}\label{sec:Ef.1}
W^{k}_\zs{r}=\{f\in \,\cC^{k}_\zs{per}[0,1]
\,,\,\sum_\zs{j=0}^k\,\|f^{(j)}\|^2\le r\}\,,
 \end{equation}
where $r>0\,,\ k\ge 1$ are some  parameters,
$\cC^{k}_\zs{per}[0,1]$ is the set of
 $k$ times continuously differentiable functions
$f\,:\,[0,1]\to\bbr$ such that $f^{(i)}(0)=f^{(i)}(1)$ for all
$0\le i \le k$. The function class $W^{k}_\zs{r}$ can be written
as an ellipsoid in $l_\zs{2}$, i.e.
 \begin{equation}\label{sec:Ef.2}
W^{k}_\zs{r}=\{f\in\,\cC^{k}_\zs{per}[0,1]\,:\,
\sum_\zs{j=1}^{\infty}\,a_\zs{j}\,\theta^2_\zs{j}\,\le r\}
 \end{equation}
where $a_\zs{j}=\sum^k_\zs{i=0}\left(2\pi [j/2]\right)^{2i}$.

Similarly to
\cite{KonevPergamenshchikov2012, KonevPergamenshchikov2015} 
we will  show here that the asymptotic sharp lower bound 
 for the robust risk \eqref{sec:In.6}
is given by
\begin{equation}\label{sec:Ef.3}
r^{*}_\zs{k}=l(r)=
\,
\left((2k+1)r\right)^{1/(2k+1)}\,
\left(
\frac{k}{(k+1)\pi} \right)^{2k/(2k+1)}\,.
\end{equation}

Note that this is  the well-known Pinsker constant
obtained for the nonadaptive filtration problem in ``signal +
small white noise'' model
 (see, for example, \cite{Pinsker1981}).

Let $\Pi_\zs{n}$ be the set of all estimators $\wh{S}_\zs{n}$
measurable with respect to the sigma-algebra
$\sigma\{y_\zs{t}\,,\,0\le t\le n\}$
 generated by the process \eqref{sec:In.1}.

\begin{theorem}\label{Th.sec:Ef.1} Under Conditions \eqref{sec:Ex.5} and \eqref{sec:Mrs.5-1}
 \begin{equation}\label{sec:Ef.4}
\liminf_\zs{n\to\infty}\,
\upsilon
^{2k/(2k+1)}_\zs{n}
 \inf_\zs{\wh{S}_\zs{n}\in\Pi_\zs{n}}\,\,
\sup_\zs{S\in W^{k}_\zs{r}} \,\cR^{*}_\zs{n}(\wh{S}_\zs{n},S) \ge
r^{*}_\zs{k}\,,
\end{equation}
where $\upsilon_\zs{n}=n/\varsigma^{*}$.
\end{theorem}
\noindent
Note that, if the parameters $r$ and $k$ are known, i.e. for the non-adaptive estimation case, in order to obtain
the efficient estimation for the ``signal+white noise'' model, Pinsker proposed in \cite{Pinsker1981} to use the estimate 
$\wh{S}_\zs{\lambda_\zs{0}}$ defined in
\eqref{sec:Mo.1} with the weights \eqref{sec:Ga.2} in which
\begin{equation}
\label{sec:Mo.11+l}
\lambda_\zs{0}=\lambda_\zs{\alpha_\zs{0}}
\quad\mbox{and}\quad
\alpha_\zs{0}=(k,\l_\zs{0})\,,
\end{equation}
where $\l_\zs{0}=[r/\varepsilon ]\varepsilon$. For the model \eqref{sec:In.1} -- \eqref{sec:Ex.1} we show the same result.

\begin{proposition}\label{Th.sec:Ef.33}
The estimator $\wh{S}_{\lambda_\zs{0}}$ satisfies the following asymptotic upper bound
$$
\lim_\zs{n \to \infty } \upsilon^{2k /(2k+1)}_\zs{n}\, \sup_\zs{S\in W^{k}_\zs{r}} \cR^*_n (\wh{S}_{\lambda_\zs{0}},S) \leq r^*_\zs{k}\,.
$$
\end{proposition}
For the adaptive estimation we user the model selection procedure \eqref{sec:Mo.9}
with the parameter $\delta$ defined  as a function of $n$ satisfying  
\begin{equation}\label{sec:Ef.4-01}
\lim_\zs{n \longrightarrow \infty }\,\delta_\zs{n}=0
\quad\mbox{and}\quad
\lim_\zs{n \longrightarrow \infty }\,n^{\check{\delta}}\,\delta_\zs{n}=0
 \end{equation}
for any $\check{\delta}>0$. For example, we can take $\delta_\zs{n}=(6+\ln n)^{-1}$.\\
\bigskip

\begin{theorem}\label{Th.sec:Ef.2}
Assume that Conditions $\H_\zs{1})$--$\H_\zs{5})$
hold true. Then the robust risk defined in  \eqref{sec:In.6} through the distribution family \eqref{sec:Ex.5}--\eqref{sec:Mrs.5-1} 
for the procedure \eqref{sec:Mo.9}  with the coefficients
\eqref{sec:Ga.2} and the parameter $\delta=\delta_\zs{n}$ satisfying \eqref{sec:Ef.4-01}
 has
 the following asymptotic upper bound
 \begin{equation}\label{sec:Ef.5}
\limsup_\zs{n\to\infty}\,
\upsilon
^{2k/(2k+1)}_\zs{n}\,
 \sup_\zs{S\in W^k_r}\,
\cR^{*}_\zs{n}(\wh{S}_\zs{*},S) \le  
r^{*}_\zs{k}
\,.
 \end{equation}
\end{theorem}
\medskip
\noindent Theorem~\ref{Th.sec:Ef.1} and Theorem~\ref{Th.sec:Ef.2} imply the following result.

\begin{corollary}\label{Co.sec:Mr.1}
Under the conditions of Theorem~\ref{Th.sec:Ef.2},
\begin{equation}\label{sec:Ef.6}
\lim_\zs{n\to\infty}\,
\upsilon
^{2k/(2k+1)}_\zs{n}\,
 \inf_\zs{\wh{S}_\zs{n}\in\Pi_\zs{n}}\,\,
\sup_\zs{S\in W^{k}_\zs{r}} \,\cR^{*}_\zs{n}(\wh{S}_\zs{n},S)
= r^{*}_\zs{k}\,.
 \end{equation}
\end{corollary}

 \begin{remark}
 \label{Re.sec.Mrs.1}
It is well known that
 the optimal (minimax) risk convergence rate
for the Sobolev ball $W^{k}_\zs{r}$
 is $n^{2k/(2k+1)}$ (see, for example, \cite{Pinsker1981}, \cite{Nussbaum1985}).
We see here that the efficient robust rate is
 $\upsilon^{2k/(2k+1)}_\zs{n}$, i.e. if  the distribution upper bound  $\varsigma^{*}\to 0$  as $n\to\infty$
 we obtain a faster rate with respect to $n^{2k/(2k+1)}$, and if $\varsigma^{*}\to \infty$  as $n\to\infty$
we obtain a slower rate. In the case when $\varsigma^{*}$ is constant the robuste rate is the same as the classical non robuste convergence rate.
\end{remark}

\section{Properties of the regression model \eqref{sec:In.1}}\label{sec:Prsm}


In order to prove the oracle inequalities we  need to study  the conditions
introduced in \cite{KonevPergamenshchikov2012} for the general semi-martingale model \eqref{sec:In.1}. 
To this end, we set for any $x\in\bbr^{n}$ the functions
\begin{equation}\label{sec:Prsm.1}
 B_\zs{1,Q}(x)=
 \sum_\zs{j=1}^{n} x_\zs{j} \, 
 \left( \E_\zs{Q}\xi^2_\zs{j,p} - \sigma_\zs{Q}\right)
 \quad\mbox{and}\quad
 B_\zs{2,Q}(x)
=
 \sum_\zs{j=1}^{n}\,x_\zs{j}\,\wt{\xi}_\zs{j,p}
\,,
\end{equation}
where $\sigma_\zs{Q}$ is defined in
\eqref{sec:Ex.5}
and $\wt{\xi}_\zs{j,p}=\xi^2_\zs{j,p}- \E_\zs{Q}\xi^2_\zs{j,p}$.

\begin{proposition}\label{Le.sec:A.06-2}
Assume that Conditions $\H_\zs{1})$--$\H_\zs{4})$ hold true. Then
\begin{equation}\label{sec:L_2-Upp_1-1}
 \L_\zs{1,Q} =
 \sup_\zs{p\geq 3}\ \sup_\zs{x\in [-1,1]^{n}}\,
  \left\vert
  B_\zs{1,Q}(x)
 \right\vert
  <2 \check{\tau}\,
 \Vert\Upsilon\Vert_\zs{1}
  \,
  \sigma_\zs{Q}
  \,.
 \end{equation}
\end{proposition}
\proof
Firstly, we set
\begin{equation}\label{sec:Pr.15}
I_n^w(f) = \int_\zs{0}^{n} f(t) \d w_\zs{t}, \, \, I_n^L(f) = \int_\zs{0}^{n} f(t) \d L_\zs{t}
\quad\mbox{and}\quad
 I_n^z(f) = \int_\zs{0}^{n} f(t) \d z_\zs{t}\,.
\end{equation}
Then,
$$
\xi_\zs{j,n} = \frac{\varrho_\zs{1}}{\sqrt n} I_n^w (\Psi_\zs{j,p})  + \frac{\varrho_\zs{2}}{\sqrt n} I_n^L (\Psi_\zs{j,p})  + \frac{\varrho_\zs{3}}{\sqrt n} I_n^z (\Psi_\zs{j,p})
$$
In view of \eqref{sec:Ex.2} the last integral can be represented as
\begin{equation}\label{sec:Pr.16}
I_n^z(f)=
 \sum_\zs{l=1}^\infty\,f(T_\zs{l}) Y_\zs{l} \Chi_\zs{\{T_\zs{l} \leq n\}}
 \,.
\end{equation}
Therefore,
\begin{equation}\label{sec:Pr.19}
\E{\xi^2_\zs{j,n}} = \frac{\bar{\varrho}}{ n} \int_\zs{0}^{n} \Psi^2_\zs{j,p}(t) \d\,t + 
\frac{\varrho^2_3}{ n} \E \sum_\zs{l=1}^\infty \Psi^2_\zs{j,p}(T_\zs{l}) 1_\zs{\{T_\zs{l} \leq n\}}\,.
\end{equation}
Using  Proposition $5.2$ from \cite{VladSlimSerge2016} we get
\begin{align*}
\E \sum_\zs{l=1}^\infty \Psi^2_\zs{j,p}(T_\zs{l}) 1_\zs{\{T_\zs{l} \leq n\}}\,&=\,
\int^{n}_\zs{0}\,\Psi^2_\zs{j,p}(x)\,\rho(x) \d\,x\\
&=\frac{1}{\check{\tau}}\,\int_\zs{0}^{n} \Psi^2_\zs{j,p}(x) \d\,x\,
+\int_\zs{0}^{n} \Psi^2_\zs{j,p}(x)\Upsilon(x) \d\,x\,,
\end{align*}
where $\rho$ is the renewal density introduced in \eqref{sec:Cns.1}. Taken into account that
\begin{equation}\label{sec:Pr.19-01}
\sup_\zs{j\ge 1}\,
\left\vert \int_\zs{0}^{n} \Psi^2_\zs{j,p}(x)\Upsilon(x) \d x
\right\vert
\le\, 2\Vert \Upsilon \Vert_\zs{1}\,,
\end{equation}
Then we obtain,
$$
\left\vert \E_\zs{Q}\xi^2_\zs{j,n}-\sigma_\zs{Q} \right\vert 
\leq \frac{2 \varrho^2_3 \Vert \Upsilon \Vert_\zs{1}}{n} 
$$
where $\sigma_\zs{Q}=\bar{\varrho}+ \varrho^2_\zs{2}/\check{\tau}$. This directly implies the desired result. 
\fdem

To study
 the function
$B_\zs{2,Q}(x),$
we have to analyze the correlation properties for  the following stochastic integrals
\begin{equation}\label{sec:tilde_Int}
 \wt{I}_\zs{n}(f)= I^2_n(f)- \E I_n^2(f)\,.
\end{equation}
To do this we 
set
\begin{equation}
\label{sec:constant_C_1}
\check{c}_\zs{1}=1+\Pi(x^{4})
+\Vert\Upsilon\Vert^{2}_\zs{1}
+\vert\rho\vert_\zs{*}
\quad\mbox{and}\quad
\check{c}_\zs{2}
=
12 (1+\check{\tau})^{2}\,(1+\check{c}_\zs{1})
\,.
\end{equation}

\noindent  
Now we investigate the behavior of the integrals defined in \eqref{sec:tilde_Int}
as functions of $f$.

\begin{proposition}\label{Le.sec:A.05-01} 
For any left continuous functions $f,g:(0,\infty) \longrightarrow \mathbb{R}$ such that 
$ \Vert f \Vert_* \leq 1$, $ \Vert g \Vert_* \leq 1$, we have 
\begin{equation}\label{sec:Upper_bound_Corls_1}
\vert \E \wt{I}_\zs{n}(f) \wt{I}_\zs{n}(g)
\vert
\le 
12 \sigma^{2}_\zs{Q}(1+\check{\tau})^{2}
\left(
(f,g)^{2}_\zs{n}
+n\check{c}_\zs{1}
\right)\,.
\end{equation}
\end{proposition}

\noindent 
Using these properties we can obtain the following bound.

\begin{proposition}\label{Le.sec:A.06-3}
Assume that Conditions $\H_\zs{1})$--$\H_\zs{4})$ hold true. Then, for all $ n\ge\,1$,
\begin{equation}\label{sec:L_2-Upp_2}
 \L_\zs{2,Q} =
  \sup_\zs{p\geq 3}\ \sup_\zs{\vert x\vert\le 1}\,
  \E
  \,
  B^{2}_\zs{2,Q}(x)
\le 
\check{c}_\zs{2}
\,
 \sigma^{2}_\zs{Q}\,,
  \end{equation}
where $|x|^2 = \sum_\zs{j=1}^{n} x^2_\zs{j}$.
\end{proposition}

\proof
Note that
$$ 
\E
\left(\sum_{j=2}^n x_j \wt{\xi}_{j,p}\right)^2 \le  
\frac{1}{n^2} \sum_{j=1}^n\sum_{l=1}^n \vert x_\zs{j}\vert \, \vert x_\zs{l}\vert \vert \E\wt{I}_\zs{n}(\Psi_{j,p}) \wt{I}_\zs{n}(\Psi_{l,p})\vert\,.
$$
Using here  Proposition \ref{Le.sec:A.05-01} and taking into account that 
$$
(\Psi_{j,p}\,,\,\Psi_{l,p})_\zs{n}=\int^{n}_\zs{0} \Psi_\zs{j,p}(t)\Psi_\zs{l,p}(t) \d t
=n\Chi_\zs{\{j=l\}}
\,,
$$
we obtain the bound \eqref{sec:L_2-Upp_2}.
Hence we obtain the desired result.
\fdem

Now  we can  study the estimate \eqref{sec:Mo.9}.

\begin{proposition}\label{Pr.sec:Si.1}
Assume that Conditions $\H_\zs{1})$ and $\H_\zs{4})$ hold true for
the model \eqref{sec:In.1} and that $S(\cdot)$ is continuously
differentiable. 
Then, for any $n\ge 2$ and $p\ge 3$,
\begin{equation}\label{sec:Si.3}
\E_\zs{Q,S}|\wh{\sigma}_\zs{n}-\sigma_\zs{Q}|
\le
\check{c}_\zs{3} \left (\frac{\sqrt{n}}{\check{p}}+\frac{1}{\sqrt{\check{p}}} \right) 
(1+\Vert\dot{S}\Vert^2)(1+\sigma_\zs{Q})^{2}
\,,
\end{equation}
where $\check{c}_\zs{3}= 6\left(14
+2\vert\rho\vert_\zs{*}
+3\sqrt{1+\check{c}_\zs{1}}\right)(1+\check{\tau})
$.
\end{proposition}

\begin{remark}\label{Re.sec:Prsm.1}
Propositions  \ref{Le.sec:A.06-2} and
\ref{Le.sec:A.06-3}
are used to obtain
the oracle inequalities given in Section \ref{sec:Mrs}
 (see, for example, \cite{KonevPergamenshchikov2012}).
\end{remark}

\section{Simulation}\label{sec:Siml}

In this section we report the results of a Monte Carlo experiment to assess the performance of the proposed model selection procedure \eqref{sec:Mo.9}. 
In \eqref{sec:In.1} we chose a 
 $1$-periodic function which, for $0\le t\le 1,$ is defined as 
\begin{equation}\label{sec:Siml.0}
S(t)=
\left\{
\begin{array}{ll}
&\vert  t-\frac{1}{2}  \vert\,
\quad\mbox{if}\quad
\,  \frac{1}{4}\leq  t \leq \frac{3}{4} \,,\\[3mm]
& \frac{1}{4}
\quad\mbox{elsewhere.}\quad
\end{array}
\right.
\end{equation}
We simulate the model
$$
\d y_\zs{t} = S(t) \d t + \d\xi_\zs{t}\,, 
$$
where $\xi_t= 0.5 w_\zs{t}+ 0.5 z_\zs{t}$. Here $z_\zs{t}$ is the semi-Markov process defined in \eqref{sec:Ex.2} with a Gaussian $\mathcal{N}(0,1)$ 
sequence $(Y_\zs{j})_\zs{j\geq1}$ and $(\tau_k)_{k\geq1}$ used in \eqref{sec:Ex.4} taken as
$\tau_k \sim \chi_\zs{3}^2$ .

We use the model selection procedure \eqref{sec:Mo.9} with the weights \eqref{sec:Ga.2}
in which $k^*= 100+\sqrt(\ln(n))$, $t_\zs{i}=i/ \ln (n)$, $m=[\ln^2 (n)]$ and $\delta=(3+\ln(n))^{-2}$.
We define the empirical risk as
\begin{equation}\label{sec:Siml.1}
\mathbf{\overline{R}}= \frac{1}{p} \sum_{j=1}^{p}  \mathbf{\hat{E}} \left(\hat{S}_n(t_j)-S(t_j)\right)^2\,,
\end{equation}
where the observation frequency   $p=100001$ and  the expectations was taken as an average over $N= 10000$ replications, i.e.
$$
\mathbf{\hat{E}} \left(\hat{S}_n(.)-S(.)\right)^2 = \frac{1}{N} \sum_{l=1}^{N} \left(\hat{S}^l_n(\cdot)-S(\cdot) \right)^2
\,.
$$
\noindent 
We set the relative quadratic risk as
\begin{equation}\label{sec:Siml.2}
\mathbf{\overline{R_*}}=\mathbf{\overline{R}}/ \Vert S\Vert ^2_\zs{p}
\quad\mbox{and}\quad
\Vert S\Vert ^2_p = \frac{1}{p} \sum_{j=0}^p S^2(t_j)\,.
\end{equation}
In our case $\Vert S\Vert ^2_p = 0.1883601$. The table below gives the values for the sample risks \eqref{sec:Siml.1} and \eqref{sec:Siml.2} for different numbers of observations $n$.
{\renewcommand{\arraystretch}{2} 
{\setlength{\tabcolsep}{1cm} 
\begin{table}
\begin{center}
    \begin{tabular}{|l|c|r|}
                                  \hline
                                  n & $\mathbf{\overline{R}}$  & $\mathbf{\overline{R_*}}$  \\
                                  \hline
                                  20 &0.0398  & 0.211 \\
                                  \hline
                                  100 & 0.0091 &0.0483  \\
                                  \hline
                                  200 & 0.0067 &0.0355  \\
                                  \hline
                                  1000 &0.0022  &0.0116  \\
                                  \hline
\end{tabular}
\end{center}
\caption{Empirical risks} \label{tab:1}
\end{table}

Figures \ref{fig1}--\ref{fig4} show  the behavior of the regression function and its estimates
by the model selection procedure \eqref{sec:Mo.9}
depending on the values of observation periods $n$. The black full line is the regression function \eqref{sec:Siml.0}
and the red dotted line is the associated  estimator.

\begin{remark}\label{Re.sec:Mnc.1}
From numerical simulations  of the procedure \eqref{sec:Mo.9}
with various observations numbers $n$ we may conclude that
 the quality of the proposed procedure  is good 
for practical needs, i.e. for reasonable (non large)  number of observations. We can also add that the quality of the estimation improves as the number of observations increases.
\end{remark}

\newpage

\begin{figure}[!h]
\begin{center}
\includegraphics[scale=.4]{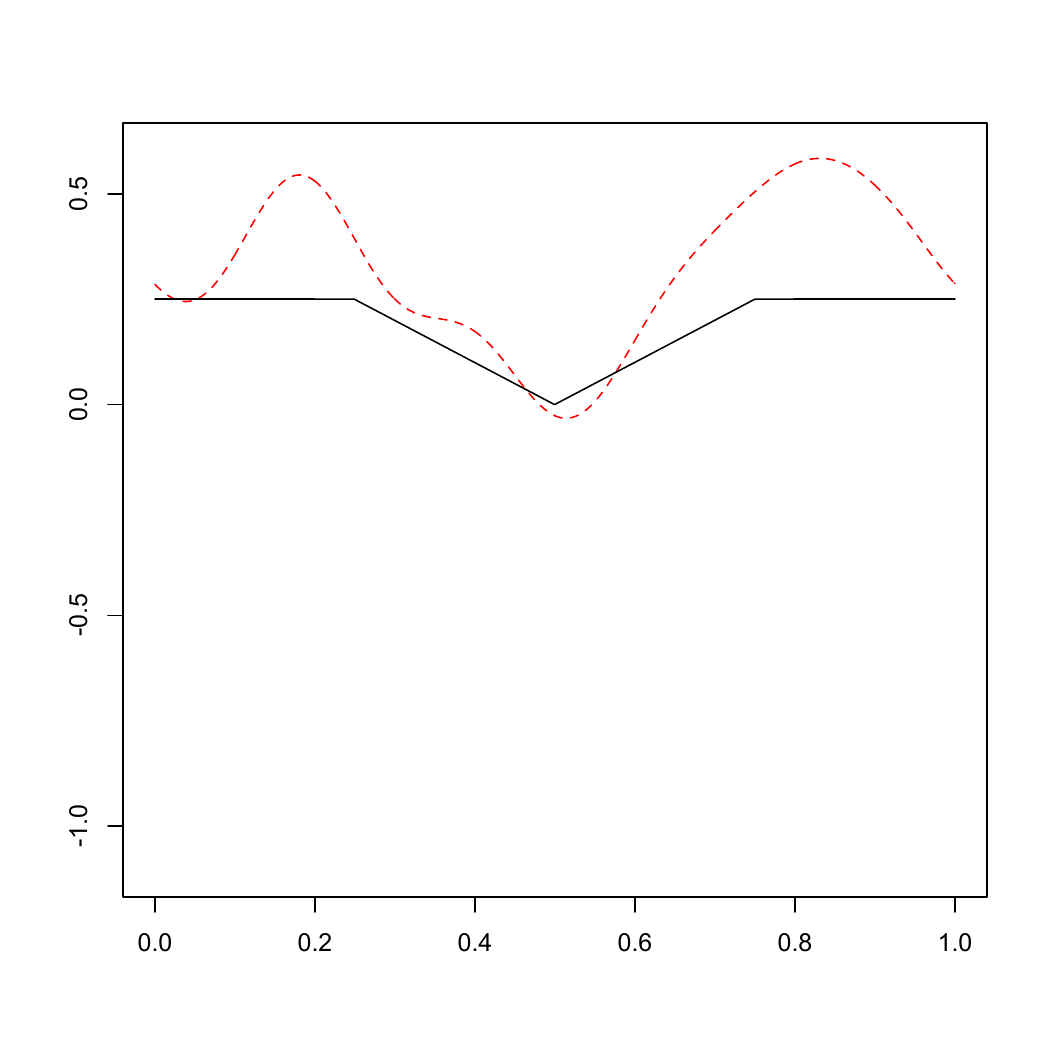}
\caption{Estimator of $S$ for $n=20$} \label{fig1}
\includegraphics[scale=0.4]{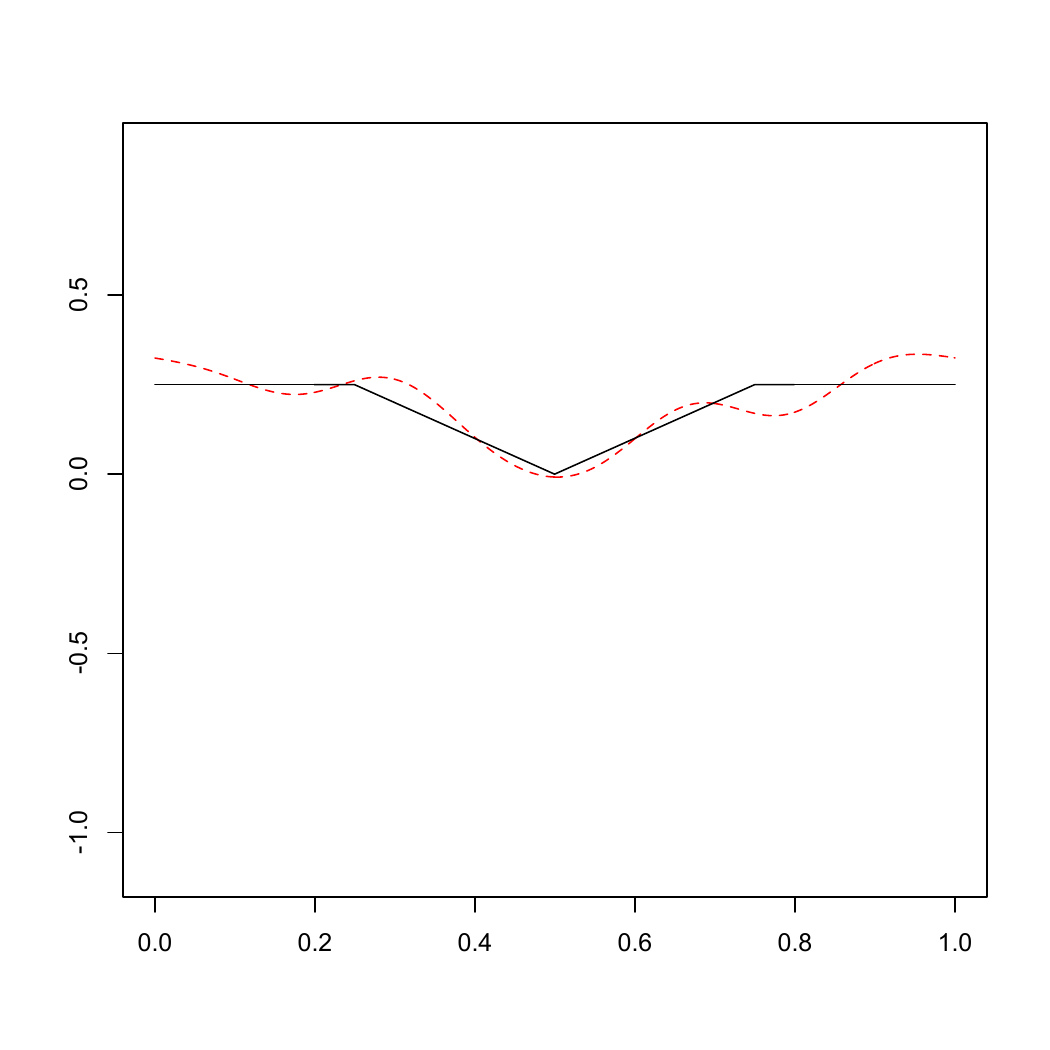} \label{fig2}
\caption{Estimator of $S$ for $n=100$} 
\end{center}
\end{figure}

\newpage

\begin{figure}[!h]
\begin{center}
\includegraphics[scale=0.4]{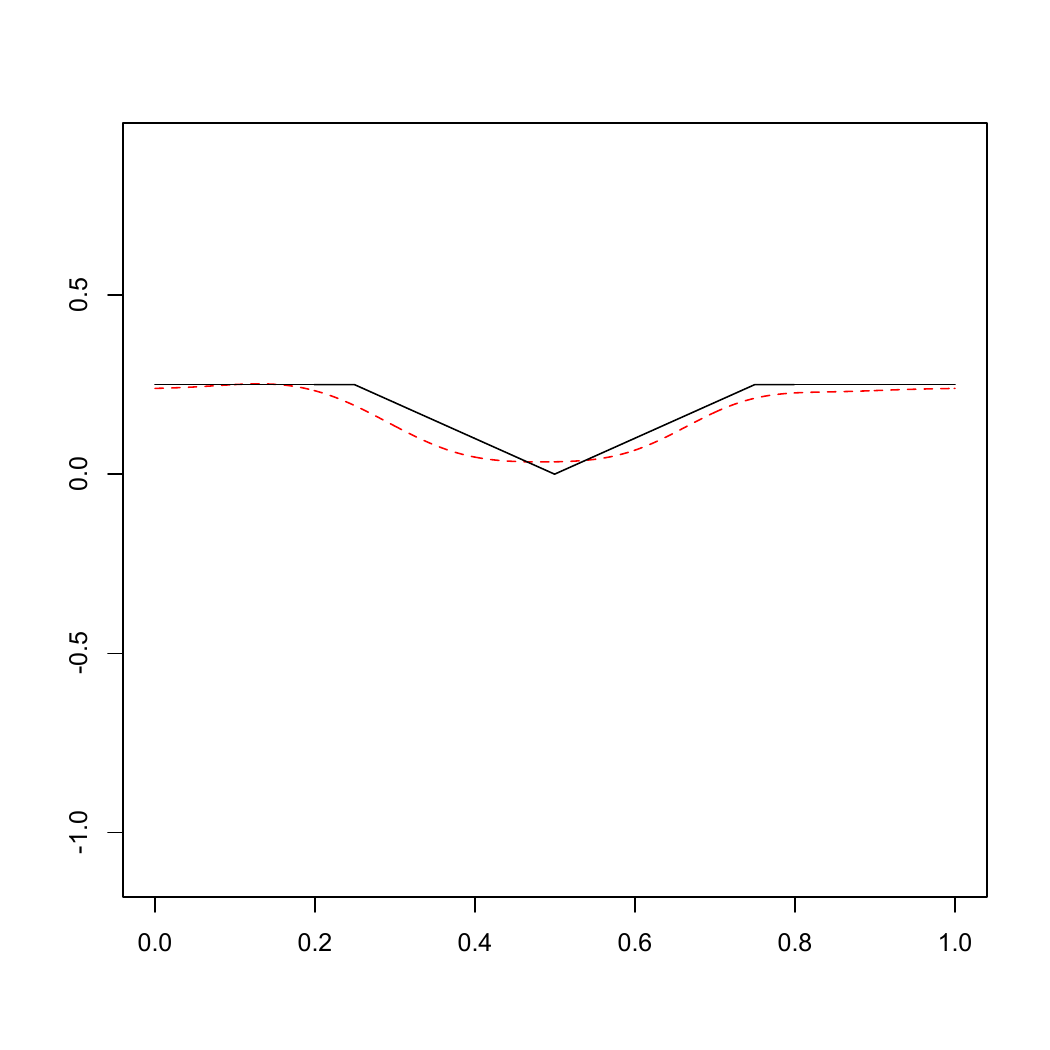} \label{fig3}
\caption{Estimator of $S$ for $n=200$} 
\includegraphics[scale=0.4]{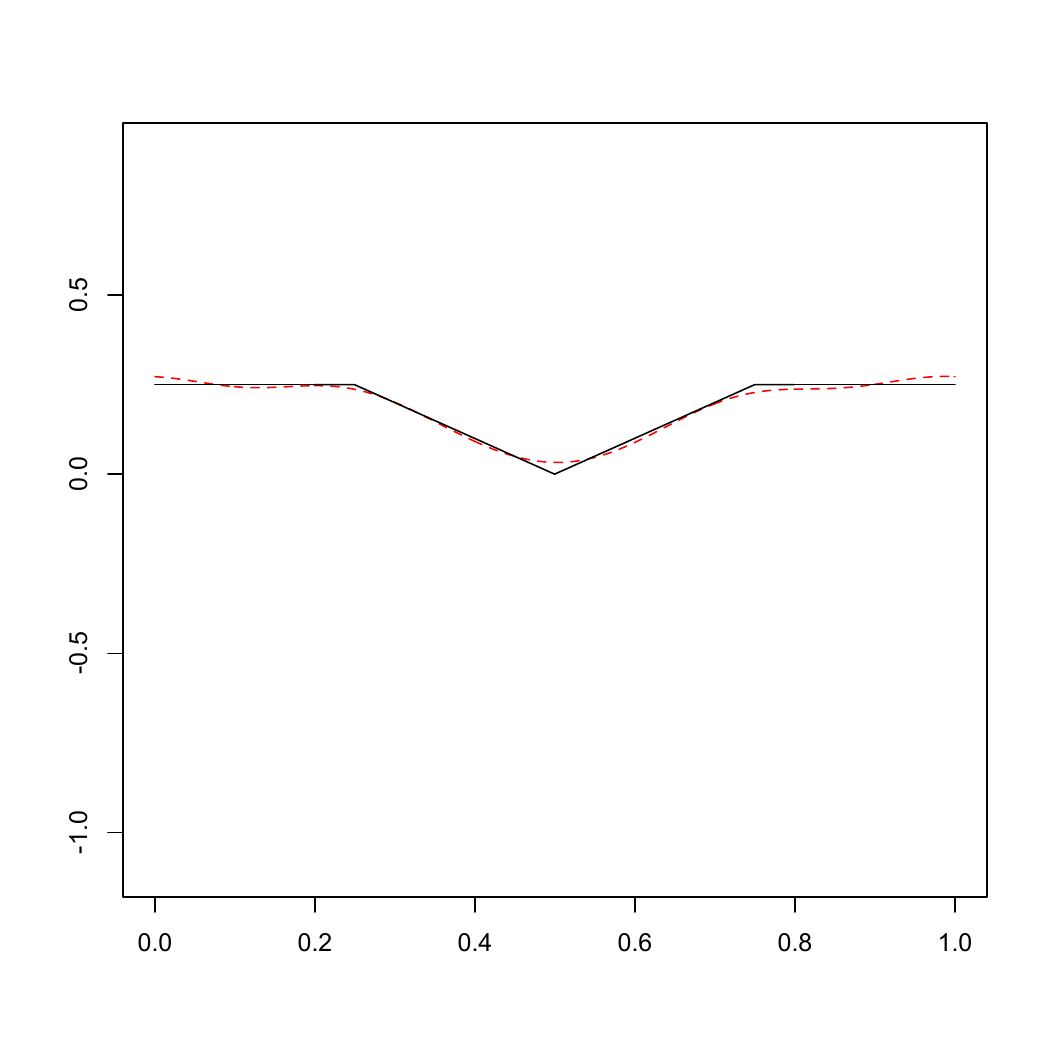} \label{fig4}
\caption{Estimator of $S$ for $n=1000$} 
\end{center}
\end{figure}
%
%
%
%


\newpage

\section{Proofs}\label{sec:Pr}

\subsection{Proof of Theorem \ref{Th.sec:Mrs.0}}

Using the cost function given in \eqref{sec:Mo.5}, we can rewrite the empirical squared error in \eqref{sec:Mo.3} as follows
\begin{equation}\label{sec:Pr.1}
\Er(\lambda) = J_n(\lambda) + 2 \sum_\zs{j=1}^{n} \lambda(j) \check{\theta}_\zs{j,p}+ \Vert S\Vert ^2-\rho \hat P_n(\lambda),
\end{equation}
where
$$\check{\theta}_\zs{j,p}=\wt{\theta}_\zs{j,p}-\overline{\theta}_\zs{j,p}\wh{\theta}_\zs{j,p}= \frac{1}{\sqrt n}\overline{\theta}_\zs{j,p}\xi_\zs{j,p} +\frac{1}{n}
\wt{\xi}_\zs{j,p}  + \frac{1}{n} \varsigma_\zs{j,n}  +\frac{\sigma_\zs{Q} -  \wh{\sigma}_\zs{n} }{n},$$
with
$\varsigma_\zs{j,p}=\E_\zs{Q}\xi^{2}_\zs{j,p}-\sigma_\zs{Q}$ and $\wt{\xi}_\zs{j,p}=\xi^{2}_\zs{j,p}-\E_\zs{Q}\xi^{2}_\zs{j,p}$. Setting
\begin{equation}\label{sec:Pr.2}
M(\lambda) = \frac{1}{\sqrt n}\sum_\zs{j=1}^{n} \lambda(j)\theta_\zs{j} \xi_\zs{j,p}
\quad\mbox{and}\quad
P^{0}_\zs{n}=\frac{\sigma_\zs{Q}\vert\lambda\vert^{2}}{n}\,,
\end{equation}
we can rewrite \eqref{sec:Pr.1} as
\begin{align}\nonumber
\Er(\lambda) & =  J_n(\lambda) + 2 \frac{\sigma_\zs{Q}-  \wh{\sigma}_\zs{n} }{n}\,L(\lambda)+ 2 M(\lambda)+\frac{2}{n}  B_\zs{1,Q}(\lambda)\\  \label{sec:Pr.3}
&+  2 \sqrt{P^{0}_n(\lambda)} \frac{ B_\zs{2,Q}(e(\lambda))}{\sqrt{\sigma_\zs{Q} n}} + \Vert S\Vert^2-\rho  P_n(\lambda),
\end{align}
where $e(\lambda)=\lambda/|\lambda|$ and the function $L(\cdot)$ was defined in \eqref{sec:Mo.2}. Let $\lambda_0= (\lambda_0(j))_\zs{1\le j\le\,p}$ be a fixed sequence in $\Lambda$ and $\wh{\lambda}$ be defined as in \eqref{sec:Mo.8}.
Substituting $\lambda_0$ and $\wh{\lambda}$ in Equation \eqref{sec:Pr.3}, we obtain
\begin{align}\label{sec:Pr.4}
\Er(\wh{\lambda})-\Er(\lambda_0)= &  J(\wh{\lambda})-J(\lambda_0)+
2 \frac{\sigma_\zs{Q}-\wh{\sigma}_\zs{n}}{n}\,L(\varpi)
+ \frac{2}{n}  B_\zs{1,Q}(\varpi)+2 M(\varpi)\nonumber\\[2mm]
& + 2 \sqrt{P^{0}_\zs{n}(\wh{\lambda})} \frac{ B_\zs{2,Q}(\wh e)}{\sqrt{\sigma_\zs{Q} n}}-2 \sqrt{P^{0}_\zs{n}(\lambda_0)}
 \frac{ B_\zs{2,Q}(e_0)}{\sqrt{\sigma_\zs{Q} n}}\nonumber \\[2mm]
& -  \delta  P_n(\wh{\lambda})+\delta P_n(\lambda_0),
\end{align}
where $\varpi= \wh{\lambda} - \lambda_\zs{0}$, $\wh{e} = e(\wh{\lambda})$ and $e_0 = e(\lambda_0)$. Note that, by \eqref{sec:Mo.2},
$$ 
|L(\varpi)| \le\,L(\hat \lambda) + L(\lambda) \leq 2\vert\Lambda\vert_\zs{*}
\,. 
$$
The inequality
\begin{equation}\label{sec:Pr.5}
2|ab| \leq \delta a^2 + \delta^{-1} b^2
\end{equation}
implies that, for any $\lambda\in\Lambda,$
$$
2 \sqrt{P^{0}_n(\lambda)} \frac{| B_\zs{2,Q}(e(\lambda))|}{\sqrt{\sigma_\zs{Q} n}} \le\, \delta P^{0}_\zs{n}(\lambda) +
 \frac{B^2_\zs{2,Q}(e(\lambda))}{\delta\sigma_\zs{Q}\,n}.
$$
Taking into account that $0 < \delta < 1$, we get
\begin{align*}
\Er(\hat \lambda) & \le\,\Er(\lambda_0) +2 M(\varpi)+ \frac{2 \L_\zs{1,Q}}{n}+ \frac{2 B^*_\zs{2,Q}}{\delta\sigma_\zs{Q}\,n} \\[2mm]
& +  \frac{1}{n} |\wh{\sigma}_\zs{n} -\sigma_\zs{Q}| ( |\wh{\lambda}|^2 + |\lambda_0|^2)+ 2  \delta P_n(\lambda_0)\,,
\end{align*}
where $B^*_\zs{2,Q} = \sup_\zs{\lambda\in\Lambda} B^2_\zs{2,Q}((e(\lambda))$. 
Moreover, noting that in view of \eqref{sec:Mo.2} $\sup_\zs{\lambda\in\Lambda} |\lambda|^2 \leq \vert\Lambda\vert_\zs{*}$,
we can rewrite the previous  bound as
\begin{align}\label{sec:Pr.6}
\Er(\wh{\lambda})  \le &  \Er(\lambda_0) +2 M(\varpi)
+ \frac{2 \L_\zs{1,Q}}{n}
+ \frac{2 B^*_\zs{2,Q}}{\delta\sigma_\zs{Q} n}  \nonumber\\[2mm]
& +  \frac{4\vert\Lambda\vert_\zs{*}}{n} |\wh{\sigma} -\sigma_\zs{Q}| + 2  \delta P_n(\lambda_0).
\end{align}
To estimate the second term in the right side of this inequality we set
$$
S_\zs{x} = \sum_\zs{j=1}^{n} x(j) \overline{\theta}_\zs{j,p} \phi_\zs{j}\,,
\quad x=(x(j))_\zs{1\le j\le n}\in\bbr^{n}\,.
$$
Thanks to 
\eqref{sec:In.3} we estimate the term $M(x)$ for any $x\in\bbr^{n}$ as
\begin{equation}\label{sec:Pr.8}
\E_\zs{Q} M^2 (x) \leq \varkappa_\zs{Q} \frac{1}{n} \sum_\zs{j=1}^{n} x^2(j) \overline{\theta}^2_\zs{j,p} = \varkappa_\zs{Q}\frac{1}{n} \Vert S_\zs{x}\Vert^2.
\end{equation}
To estimate this function for a random vector $x\in\bbr^{n}$, we set
$$ 
Z^* = \sup_\zs{x \varepsilon \Lambda_1} \frac{n M^2 (x)}{\Vert S_x\Vert^2}\,,
\quad
\Lambda_1 = \Lambda - \lambda_0\,.
$$
So, through the Inequality  \eqref{sec:Pr.5}, we get 
\begin{equation}\label{sec:Pr.10}
2 |M(x)|\leq \delta \Vert S_x\Vert^2 + \frac{Z^*}{n\delta}.
\end{equation}
It is clear that the last term  here can be estimated as
\begin{equation}\label{sec:Pr.9}
\E_\zs{Q} Z^* \leq \sum_\zs{x \in \Lambda_1} \frac{n \E_\zs{Q} M^2 (x)}{\Vert S_x\Vert^2} \leq \sum_\zs{x \in \Lambda_1} \varkappa_\zs{Q}= \varkappa_\zs{Q}\nu\,,
\end{equation}
where $\nu = \mbox{card}(\Lambda)$. 
Moreover, note that, for any $x\in\Lambda_1$,
\begin{equation}\label{sec:Pr.11}
\Vert S_x\Vert^2-\Vert\wh{S}_x\Vert^2 = \sum_\zs{j=1}^{n} x^2(j) (\overline{\theta}^2_\zs{j,p}-\wh{\theta}^2_\zs{j,p}) \le -2 M_1(x),
\end{equation}
where $M_\zs{1}(x) =  n^{-1/2}\,\sum_\zs{j=1}^{n}\, x^2(j)\overline{\theta}^2_\zs{j,p} \xi_\zs{j,n}$.
Taking into account now that, for any $x \in \Lambda_1$, the components $|x(j)|\leq 1$,  we can estimate this term as
in \eqref{sec:Pr.8}, i.e.
$$
\E_\zs{Q}\, M^2_\zs{1}(x) \leq \varkappa_\zs{Q}\,
\frac{\Vert S_x\Vert^2}{n}\,.
$$
Similarly to the previous reasoning
we set
$$ 
Z^*_\zs{1} = \sup_\zs{x \varepsilon \Lambda_1} \frac{n M^2_1 (x)}{\Vert S_x\Vert^2}
$$
and we get
\begin{equation}\label{sec:Pr.12}
\E_\zs{Q}\, Z^*_1 \leq \varkappa_\zs{Q}\,\nu\,.
\end{equation}
Using the same type of arguments as in \eqref{sec:Pr.10}, we can derive
\begin{equation}\label{sec:Pr.13}
2 |M_1(x)|\leq \delta \Vert S_x\Vert^2 + \frac{Z^*_1}{n\delta}.
\end{equation}
From here and \eqref{sec:Pr.11}, we get
\begin{equation}\label{sec:Pr.14}
\Vert S_x\Vert^2 \leq \frac{\Vert\wh{S}_x\Vert^2}{1-\delta} + \frac{Z^*_1}{n \delta (1-\delta)}
\end{equation}
for any $0<\delta<1$. Using this bound in \eqref{sec:Pr.10} yields
$$
2 M(x) \leq \frac{\delta \Vert\wh{S}_x\Vert^2}{1-\delta} + \frac{Z^*+Z^*_1}{n \delta (1-\delta)}
\,.
$$
Taking into account that $\Vert\wh{S}_\zs{\varpi}\Vert^{2}\le 2\,(\Er(\wh{\lambda})+\Er(\lambda_0))$, we obtain
$$
2 M(\varpi) \leq \frac{2\delta(\Er(\wh{\lambda})+\Er(\lambda_0))}{1-\delta} + \frac{Z^*+Z^*_1}{n \delta (1-\delta)}.
$$
Using this bound in  \eqref{sec:Pr.6} 
we obtain 
\begin{align*}
\Er(\wh{\lambda}) \le &  \frac{1+\delta}{1-3\delta} \Er(\lambda_0) 
+ \frac{Z^*+Z^*_1}{n \delta (1-3\delta)}
+ \frac{2 \L_\zs{1,Q}}{n(1-3\delta)}
+ \frac{2 B^*_\zs{2,Q}}{\delta(1-3\delta)\sigma_\zs{Q} n} \\[2mm]
& +  \frac{(4\vert\Lambda\vert_\zs{*}+2)}{n(1-3\delta)} |\wh{\sigma} -\sigma_\zs{Q}| + \frac{2\delta}{(1-3\delta)} P^{0}_n(\lambda_0).
\end{align*}
Moreover, for $0<\delta<1/6$ we can rewrite this inequality as
\begin{align*}
\Er(\wh{\lambda}) \le &  \frac{1+\delta}{1-3\delta} \Er(\lambda_0) 
+ \frac{2(Z^*+Z^*_1)}{n \delta}
+ \frac{4 \L_\zs{1,Q}}{n}
+ \frac{4 B^*_\zs{2,Q}}{\delta \sigma_\zs{Q} n} \\[2mm]
& +  \frac{(8\vert\Lambda\vert_\zs{*}+2)}{n} |\wh{\sigma} -\sigma_\zs{Q}| + \frac{2\delta}{(1-3\delta)}\,
 P^{0}_n(\lambda_0).
\end{align*}
Now, in view of the condition Proposition~\ref{Le.sec:A.06-3}, we estimate the expectation of the term $B^*_\zs{2,Q}$ in \eqref{sec:Pr.6} as
$$
\E_\zs{Q}\, B^*_\zs{2,Q} \leq \sum_\zs{\lambda\in\Lambda}\E_\zs{Q} B^2_\zs{2,Q} (e(\lambda)) \leq \nu \L_\zs{2,Q}\,.
$$
Now,  taking into account that $\vert\Lambda\vert_\zs{*}\ge 1$, we get
\begin{align*}
 \cR_\zs{Q}(\wh{S}_*,S) \le & \frac{1+\delta}{1-3\delta} 
 \cR_\zs{Q}(\wh{S}_\zs{\lambda_0},S)
+ \frac{4\varkappa_\zs{Q} \nu}{n \delta}
+ \frac{4 \L_\zs{1,Q}}{n}
+ \frac{4 \nu \L_\zs{2,Q}}{\delta \sigma_\zs{Q} n} \\[2mm]
& +  \frac{10\vert\Lambda\vert_\zs{*}}{n} \,\E_\zs{Q}\,|\wh{\sigma} -\sigma_\zs{Q}| + \frac{2\delta}{(1-3\delta)} P^{0}_n(\lambda_0).
\end{align*}
By using the upper bound for $ P^0_n(\lambda_0)$ in Lemma~\ref{Le.sec:A.1-06-11-01}, we obtain that
\begin{align*}
 \cR_\zs{Q}(\wh{S}_*,S) \le &  \frac{1+3\delta}{1-3\delta} 
 \cR_\zs{Q}(\wh{S}_\zs{\lambda_0},S)
+ \frac{4\varkappa_\zs{Q} \nu}{n \delta}
+ \frac{4 \L_\zs{1,Q}}{n}
+ \frac{4 \nu \L_\zs{2,Q}}{\delta \sigma_\zs{Q} n} \\[2mm]
& +  \frac{10\vert\Lambda\vert_\zs{*}}{n} \,\E_\zs{Q}\,|\wh{\sigma} -\sigma_\zs{Q}| + \frac{2\delta \L_\zs{1,Q}}{(1-3\delta)n}
\,.
\end{align*}
Taking into account here that  $1-3\delta\ge 1/2$ for $0<\delta<1/3$ 
and that $\varkappa_\zs{Q}\le (1+\check{\tau}\vert\rho\vert_\zs{*})\sigma_\zs{Q}$
and using the bounds
\eqref{sec:L_2-Upp_1-1} and
\eqref{sec:L_2-Upp_2} we obtain the inequality \eqref{sec:Mrs.1--00}. Hence Theorem \ref{Th.sec:Mrs.0}\,.

 \fdem
 
\subsection{Proof of Proposition \ref{Le.sec:A.05-01} }
By Ito's formula  one gets
\begin{equation}
\label{sec:Pr.4-11-00}
\d I^2_\zs{t}(f)= 2 I_\zs{t-}(f) \d I_\zs{t}(f)+ \varrho_\zs{1}^2 f^2(t) \d\,t + \sum_\zs{0\leq s \leq t}f^2(s) (\Delta \xi^{d}_\zs{s})^2\,,
\end{equation}
where $\xi^{d}_\zs{t}=\varrho_\zs{2}\,L_\zs{t}+\varrho_\zs{3}z_\zs{t}$. Taking into account that the processes 
$(L_\zs{t})_\zs{t\ge 0}$
and $(z_\zs{t})_\zs{t\ge 0}$ are independent and the time of jumps $T_\zs{k}$ defined in \eqref{sec:Ex.4}
have a density,
we have $\Delta z_\zs{s} \Delta L_\zs{s}=0$ a.s. for any $s\ge 0$.
Therefore, we can rewrite the differential 
\eqref{sec:Pr.4-11-00} as
\begin{align}\nonumber
\d I^2_\zs{t}(f)=& 2 I_\zs{t-}(f) \d I_\zs{t}(f)+ \varrho_\zs{1}^2\, f^2(t) \d\,t 
\\[2mm] \label{sec:Pr.4-11-00++1}
&+ 
\varrho^{2}_\zs{2}\d \sum_\zs{0\leq s \leq t}f^2(s) (\Delta L_\zs{s})^2
+
\varrho^{2}_\zs{3}\d \sum_\zs{0\leq s \leq t}f^2(s) (\Delta z_\zs{s})^2
\,.
\end{align}
Therefore, using Lemma \ref{Le.sec:Smp.2}  we obtain
$$
\E I^2_\zs{t}(f) = \bar{\varrho}\,
\Vert f\Vert^2_\zs{t} + \varrho_\zs{3}^2
\Vert f\sqrt{\rho}\Vert^2_\zs{t}\,,
$$
where $\Vert f\Vert^{2}_\zs{t}=\int^{t}_\zs{0}\,f^{2}(t)\d t$, $\rho$ is the density of the renewal  measure $\sum_\zs{j=1}^{\infty} \eta^{(j)}$ and with $\eta$ the distribution
of $\tau_\zs{1}$.
Therefore,
\begin{equation}\label{sec:A.Int_12++}
\d\wt{I}_\zs{t}(f)= 2 I_\zs{t-}(f) f(t) d\xi_\zs{t}
+
f^2(t) \d \wt{m}_\zs{t}\,,\quad 
\wt{m}_\zs{t}=\varrho^{2}_\zs{2}\check{m}_\zs{t}+\varrho^{2}_\zs{3}m_\zs{t}\,,
\end{equation}
where 
$
\check{m}_\zs{t}= \sum_\zs{0\leq s \leq t}(\Delta L_\zs{s})^2 - t$ and
$
m_\zs{t} = \sum_\zs{0\leq s \leq t}(\Delta z_\zs{s})^2 -
 \int_\zs{0}^{t} \rho(s) \d s$.
 By the Ito formula we get
\begin{align}\nonumber
\E \wt{I}_\zs{n}(f) \wt{I}_\zs{n}(g) = & \E \int^n_0 \wt{I}_\zs{t-}(f) \d\wt{I}_\zs{t}(g) 
\\[2mm] \label{sec:AppInt_1}
&+ \E \int^n_0 \wt{I}_\zs{t-}(g) \d\wt{I}_\zs{t}(f)
+ \E \left[ \wt{I} (f),\wt{I} (g)  \right]_\zs{n}\,.
\end{align}
First, note that the process $(\check{m}_\zs{t})_\zs{t\ge 0}$ is a martingale and, using
Lemma \ref{Le.sec:A.05-00}, we get
$$
\E  \int^n_0 \wt{I}_\zs{t-}(f) \d\wt{I}_\zs{t}(g) =\rho_\zs{2}^2\E \int^n_0 \wt{I}_\zs{t-}(f)  g^2(t) \d m_\zs{t}
=
\rho_\zs{2}^2\E \int^n_0 \,I^{2}_\zs{t-}(f)  g^2(t) \d m_\zs{t}\,.
$$
The last integral can be represented as
$$
 \E\int^n_0 \,I^{2}_\zs{t-}(f)  g^2(t) \d m_\zs{t}
=
J_\zs{1}-
J_\zs{2}
\,,
$$
where
$$
J_\zs{1}=\E\sum_\zs{k\ge 1}\,I^{2}_\zs{T_\zs{k}-}(f)  g^2(T_\zs{k})\Chi_\zs{\{T_\zs{k}\le n\}}
\quad\mbox{and}\quad
J_\zs{2}=\int^n_0 \,\E\,I^{2}_\zs{t}(f)  g^2(t)\rho(t) \d t\,.
$$
By Lemma \ref{Le.sec:A.06-11-03} we get
$$
J_\zs{1}=\E\sum_\zs{k\ge 1}\,\E\left(I^{2}_\zs{T_\zs{k}-}(f)\vert\cG\right)  g^2(T_\zs{k})\Chi_\zs{\{T_\zs{k}\le n\}}
=
\varrho_\zs{1}^2 J_\zs{1,1}
+ 
\varrho_\zs{2}^2J_\zs{1,2}\,,
$$
where
$$
J_\zs{1,1}=\E\sum_\zs{k\ge 1}\,
\Vert f\Vert^{2}_\zs{T_\zs{k}}
g^2(T_\zs{k})
\Chi_\zs{\{T_\zs{k}\le n\}}
\quad\mbox{and}\quad
J_\zs{1,2}=
\E\sum_\zs{k\ge 1}\,
\sum_\zs{l=1}^{k-1}\, f^{2}(T_\zs{l})
\,  g^2(T_\zs{k})
\Chi_\zs{\{T_\zs{k}\le n\}}\,.
$$
We obtain directly that
$$
J_\zs{1,1}=\int^{n}_\zs{0}\,\Vert f\Vert^{2}_\zs{t}
g^2(t)\rho(t)\d t
$$
and
$$
J_\zs{1,2}=
\E
\sum_\zs{l\ge 1}\, f^{2}(T_\zs{l})
\sum_\zs{k\ge l+1}\,
\,  g^2(T_\zs{k})
\Chi_\zs{\{T_\zs{k}\le n\}}
=\int^{n}_\zs{0}\,f^{2}(x)
\left(\int^{n-x}_\zs{0}\,g^{2}(x+t)\rho(t)\d t \right)
\rho(x)\d x.
$$
From Lemma \ref{Le.sec:Smp.2} we obtain that
$$
J_\zs{2}=\bar{\varrho}\int^{n}_\zs{0}\Vert f\Vert^{2}_\zs{t}g^{2}(t)\rho(t)\d t
+
\varrho^{2}_\zs{3}\int^{n}_\zs{0}\Vert f \sqrt{\rho}\Vert^{2}_\zs{t}g^{2}(t)\rho(t)\d t
\,.
$$
Therefore,
$$
 \E\int^n_0 \,I^{2}_\zs{t-}(f)  g^2(t) \d m_\zs{t}
=\varrho^{2}_\zs{3}\int^{n}_\zs{0}f^{2}(x)
\left(
\int^{n}_\zs{x} g^{2}(t)(\rho(t-x)-\rho(t))\d t
\right)
\rho(x)\d x
\,.
$$
Taking into account that $\rho(t-x)-\rho(t)=\Upsilon(t-x)-\Upsilon(t)$
we can estimate the last integral as
$$
\vert \E\int^n_0 \,I^{2}_\zs{t-}(f)  g^2(t) \d m_\zs{t}\vert
\le 2 \varrho^{2}_\zs{3}n\Vert \Upsilon\Vert_\zs{1}\,. 
$$
From this and by the symmetry arguments we obtain that
\begin{equation}\label{sec:A_IntBound_121}
\vert \E  \int^n_0 \wt{I}_\zs{t-}(f) \d\wt{I}_\zs{t}(g)\vert
+
\vert \E  \int^n_0 \wt{I}_\zs{t-}(g) \d\wt{I}_\zs{t}(f)\vert
\le 4 \varrho^{4}_\zs{3}n\Vert \Upsilon\Vert_\zs{1}\,. 
\end{equation}

Note now that
\begin{equation}\label{sec:Pr.4-11-02}
 \left[ \wt{I} (f),\wt{I} (g)  \right]_\zs{n} = \left\langle  \wt{I}^c (f),\wt{I}^c (g)  \right\rangle_\zs{n}  + 
\D_\zs{n}(f,g)\,,
\end{equation}
where
$$
\D_\zs{n}(f,g)= \sum_{0\leq t \leq n}   \Delta \wt{I}^{d}_\zs{t}(f) \Delta \wt{I}^{d}_\zs{t}(g)\,.
$$
It should be noted that
 the continuous and the discrete parts of the processes \eqref{sec:A.Int_12++} can be represented as
$$
\wt{I}^c_\zs{t} (f)=2\varrho_\zs{1}  \int^{t}_\zs{0}\,I_\zs{s}(f) f(s) \d w_\zs{s}
\quad\mbox{and}\quad
\wt{I}^{d}_\zs{t} (g)=2 \int^{t}_\zs{0}\,I_\zs{s-}(f) f(s) d\xi^{d}_\zs{s}
+
\int^{t}_\zs{0}\,f^2(s) \d \wt{m}_\zs{s}\,.
$$
So, in view of Lemma 6.1 from \cite{VladSlimSerge2016},
\begin{align}\nonumber
\E  <  \wt{I}^c (f),&\wt{I}^c (g)  >_\zs{n}= 4 \rho_1^2  \int^{n}_0 
\E ( I_\zs{t}(f) I_\zs{t}(g))   f(t) g(t) \d t \\ \nonumber
& = 4 \rho_1^4\int^n_0  (f,g)_\zs{t}\,f(t) g(t) \d t + 
4 \rho_1^2 \rho_3^2  \int^n_0 (f,g\rho)_\zs{t} f(t) g(t) \d t\\ \label{sec:A_Int_f_g}
& =
4\rho_1^2 \sigma_\zs{Q}\,
(f,g)^{2}_\zs{n}
+
4 \rho_1^2 \rho_3^2\int^n_0 (f,g\Upsilon)_\zs{t} f(t) g(t) \d t
\,,
\end{align}
with  $(f,g)_\zs{t} = \int^t_0   f(s) g(s) \d s$. 
Taking into account that
$\Vert f \Vert_* \leq 1$ and $ \Vert g \Vert_* \leq 1,$ we can estimate the last integral as
$$
\int^n_0 (f,g\Upsilon)_\zs{t} f(t) g(t) \d t\le n \Vert\Upsilon\Vert_\zs{1}\,.
$$
Therefore, 
\begin{equation}
\label{sec:<II>__11}
\left\vert \E  \left\langle  \wt{I}^c (f),\wt{I}^c (g)  \right\rangle_\zs{n} 
\right\vert
   \leq 
4
\sigma^{2}_\zs{Q}
\left(
(f,g)^{2}_\zs{n}
+
n \check{\tau} \Vert\Upsilon\Vert_\zs{1}
 \right)  
\,.
\end{equation}
To study the last  term in \eqref{sec:Pr.4-11-02} note that
$$
\D_\zs{n}(f,g)=\sum_\zs{0\le t\le n}\left(2 I_\zs{t-}(f) f(t)\Delta \xi^{d}_\zs{t} + f^2(t) \Delta \wt{m}_\zs{t}\right)
\left(2 I_\zs{t-}(g) g(t)\Delta \xi^{d}_\zs{t} + g^2(t) \Delta \wt{m}_\zs{t}\right)
\,.
$$
Taking into account that for any $t>0$
$$
\Delta\xi^{d}_\zs{t}
\Delta\wt{m}_\zs{t}
=
\varrho^{3}_\zs{2}
(\Delta L_\zs{t})^{3}
+
\varrho^{3}_\zs{3}
(\Delta z_\zs{t})^{3}
\,,
$$
we obtain that
$$
\E\,\sum_\zs{0\le t\le n}\,I_\zs{t-}(f) f(t) g^2(t) \Delta \xi^{d}_\zs{t}\,\Delta \wt{m}_\zs{t}
=0\,.
$$
So, using the symmetry arguments, we find that
\begin{equation}\label{sec:App_D_11++}
\E\D_\zs{n}(f,g)
=4\E\,\D_\zs{1,n}(f,g)
+ \E\,
\D_\zs{2,n}(f,g)
\,,
\end{equation}
where
$$
\D_\zs{1,n}(f,g)=\sum_\zs{0\le t\le n}\, I_\zs{t-}(f)I_\zs{t-}(g) f(t)g(t)(\Delta \xi^{d}_\zs{t})^{2}
\quad\mbox{and}\quad
\D_\zs{2,n}(f,g)=\sum_\zs{0\le t\le n}\,f^2(t)\,g^2(t) (\Delta \wt{m}_\zs{t})^{2}\,.
$$
Note that 
$$
\D_\zs{1,n}(f,g)=\varrho^{2}_\zs{2}\check{\D}_\zs{1,n}(f,g)
+
\varrho^{2}_\zs{3}\wt{\D}_\zs{1,n}(f,g)
\,,
$$
where
$$
\check{\D}_\zs{1,n}(f,g)=\sum_\zs{0\le t\le n}\, I_\zs{t-}(f)I_\zs{t-}(g) f(t)g(t)(\Delta L_\zs{t})^{2}
$$
and
$$
\wt{\D}_\zs{1,n}(f,g)=\sum_\zs{0\le t\le n}\, I_\zs{t-}(f)I_\zs{t-}(g) f(t)g(t)(\Delta z_\zs{t})^{2}
\,.
$$
Now, similarly to \eqref{sec:A_Int_f_g} and taking into account that $\Pi(x^{2})=1,$ we get
\begin{align*}
\E\,\check{\D}_\zs{1,n}(f,g)=&\int^{n}_\zs{0}\,f(t)g(t)\E\, I_\zs{t}(f)I_\zs{t}(g) \,\d t=
\bar{\varrho}\,\int^{n}_\zs{0}\,f(t)g(t)\,(f,g)_\zs{t} \,\d t\\[2mm]
&+\varrho^{2}_\zs{3}
\,\int^{n}_\zs{0}\,f(t)g(t)\,(f,g\rho)_\zs{t} \,\d t\\[2mm]
=&\sigma_\zs{Q}(f,g)^{2}_\zs{n}
+\varrho^{2}_\zs{3}
\,\int^{n}_\zs{0}\,f(t)g(t)\,(f,g\Upsilon)_\zs{t} \,\d t\,.
\end{align*}
So,
\begin{equation}\label{sec:A_Upper_bound_1_1}
\vert \E\,\check{\D}_\zs{1,n}(f,g)\vert\le \sigma_\zs{Q}
\left(
(f,g)^{2}_\zs{n}
+n\check{\tau}\Vert\Upsilon\Vert_\zs{1}
\right)\,.
\end{equation}
Moreover, taking into account that $\E Y^{2}_\zs{j}=1$ we get
$$
\E\,\wt{\D}_\zs{1,n}(f,g)=\E
\sum_\zs{k\ge 1}\, I_\zs{T_\zs{k}-}(f)I_\zs{T_\zs{k}-}(g) f(T_\zs{k})g(T_\zs{k})
\,\Chi_\zs{\{T_\zs{k}\le n\}}
\,.
$$
So, in view of Lemma \ref{Le.sec:A.06-11-03} 
\begin{align*}
\E\,&\wt{\D}_\zs{1,n}(f,g)=\E
\sum_\zs{k\ge 1}\, 
\E\left( I_\zs{T_\zs{k}-}(f)I_\zs{T_\zs{k}-}(g) \vert\cG\right)
f(T_\zs{k})g(T_\zs{k})
\,\Chi_\zs{\{T_\zs{k}\le n\}}\\[2mm]
&
=\bar{\varrho}\E 
\sum_\zs{k\ge 1}(f\,,\,g)_\zs{T_\zs{k}}f(T_\zs{k})g(T_\zs{k})
\,\Chi_\zs{\{T_\zs{k}\le n\}}
+ 
\varrho_\zs{3}^2
\,\E\,\D^{'}_\zs{1,n}(f,g)
\\
&=
\bar{\varrho}
\int^{n}_\zs{0}(f,g)_\zs{t}\,f(t)g(t)\rho(t)\d t
+
\varrho_\zs{3}^2\E\,\D^{'}_\zs{1,n}(f,g)
\,,
\end{align*}
where
$$
\D^{'}_\zs{1,n}(f,g)=
\sum_\zs{k\ge 1}
 \sum_\zs{l=1}^{k-1}\, f(T_\zs{l})\,g(T_\zs{l})
 f(T_\zs{k})g(T_\zs{k})
\,\Chi_\zs{\{T_\zs{k}\le n\}}
\,.
$$
Noting now that
$$
\int^{n}_\zs{0}(f,g)_\zs{t}\,f(t)g(t)\rho(t)\d t=\frac{1}{2\check{\tau}}(f,g)^{2}_\zs{n}
+
\int^{n}_\zs{0}(f,g)_\zs{t}\,f(t)g(t)\Upsilon(t)\d t
\,,
$$
we obtain
$$
\vert\int^{n}_\zs{0}(f,g)_\zs{t}\,f(t)g(t)\rho(t)\d t\vert\le \frac{1}{2\check{\tau}}(f,g)^{2}_\zs{n}
+n\Vert\Upsilon\Vert_\zs{1}\,.
$$
Furthermore, the expectation of $\D^{'}_\zs{1,n}(f,g)$
can be represented as
\begin{align*}
\E\,\D^{'}_\zs{1,n}(f,g)&=
\E\,\sum_\zs{l\ge 1}\, f(T_\zs{l})\,g(T_\zs{l})
\sum_\zs{k\ge l+1}
 f(T_\zs{k})g(T_\zs{k})
\,\Chi_\zs{\{T_\zs{k}\le n\}}\\[2mm]
&=\int^{n}_\zs{0}f(x)g(x)\,
\left(\int^{n-x}_\zs{0}\,f(x+t)g(x+t)\rho(t)\d t
\right)
\rho(x)\d x\\[2mm]
&=
\frac{1}{2\check{\tau}}(f,g)^{2}_\zs{n}
+\D^{''}_\zs{1,n}(f,g)
\,,
\end{align*}
where the last term in this equality can be represented as
\begin{align*}
\D^{''}_\zs{1,n}(f,g)=&
\int^{n}_\zs{0}f(x)g(x)\,
\left(\int^{n-x}_\zs{0}\,f(x+t)g(x+t)\Upsilon(t)\d t
\right)
\rho(x)\d x\\[2mm]
&+
\frac{1}{\check{\tau}}\,
\int^{n}_\zs{0}f(x)g(x)\,
\left(\int^{n-x}_\zs{0}\,f(x+t)g(x+t)\Upsilon(t)\d t
\right)
\Upsilon(x)\d x
\,.
\end{align*}
This implies
$$
\vert\D^{''}_\zs{1,n}(f,g)\vert
\le n(1+\frac{1}{\check{\tau}})
(1+\Vert\Upsilon\Vert^{2}_\zs{1})
\,.
$$
Therefore,
\begin{equation}\label{sec:A_Upper_bound_1_1++}
\vert \E\,\wt{\D}_\zs{1,n}(f,g)\vert\le \sigma_\zs{Q}
\left(
(f,g)^{2}_\zs{n}
+n(1+\check{\tau})\Vert\Upsilon\Vert^{2}_\zs{1}
\right)\,.
\end{equation}
Finally we  obtain that
\begin{equation}\label{sec:A_Upper_bound_D_1}
\vert \E\,\D_\zs{1,n}(f,g)\vert\le \sigma^{2}_\zs{Q}(1+\check{\tau})^{2}
\left(
(f,g)^{2}_\zs{n}
+n\Vert\Upsilon\Vert^{2}_\zs{1}
\right)\,.
\end{equation}
As to the last term in \eqref{sec:App_D_11++} 
we can calculate directly
$$
\E\,\D_\zs{2,n}(f,g)
=\varrho^{4}_\zs{2}
\Pi(x^{4})\int^{n}_\zs{0}f^{2}(t)\,g^{2}(t)\d t
+\varrho^{4}_\zs{3}
\int^{n}_\zs{0}f^{2}(t)\,g^{2}(t)\rho(t)\d t
\,,
$$
i.e.
$$
\E\,\D_\zs{2,n}(f,g)\le n \sigma^{2}_\zs{Q}\left(\Pi(x^{4})+\vert\rho\vert_\zs{*} \right)(1+\check{\tau})^{2}
\,.
$$
From here we obtain that
\begin{equation}\label{sec:A_Upper_bound_D_n+}
\vert \E\,\D_\zs{n}(f,g)\vert\le \sigma^{2}_\zs{Q}(1+\check{\tau})^{2}
\left(
4(f,g)^{2}_\zs{n}
+n\check{c}_\zs{1}
\right)\,,
\end{equation}
where
$\check{c}_\zs{1}$ is given in
\eqref{sec:constant_C_1}.
 From this and
\eqref{sec:<II>__11}
we find
\begin{equation}\label{sec:A_Upper_bound_[]_n}
\E
[ \wt{I} (f),\wt{I} (g)]_\zs{n} 
\le 
8\sigma^{2}_\zs{Q}(1+\check{\tau})^{2}
\left(
(f,g)^{2}_\zs{n}
+n\check{c}_\zs{1}
\right)\,.
\end{equation}
This bound and \eqref{sec:A_IntBound_121} implies \eqref{sec:Upper_bound_Corls_1}.
Hence Lemma \ref{Le.sec:A.05-01}. \fdem

\subsection{Proof of Proposition \ref{Pr.sec:Si.1}}
It is clear that the Inequality \eqref{sec:Si.3} holds true for $l > \check{p}$. 
Let now $l \leqslant \check{p}$. Setting $x_j^{'}= \Chi_\zs{\{[\sqrt{n}] \leqslant  j\leqslant \check{p}\}}$ and subtituting \eqref{sec:In.8}   in \eqref{sec:Mo.4-1-31-3} yields,

\begin{equation}\label{sec:Mo.1-0-1-04}
\wh{\sigma}_\zs{n}=
\frac{n}{\check{p}} \sum^{\check{p}}_\zs{j=l}\, (\overline{\theta}_\zs{j,p})^2
+
\frac{2n}{\check{p}} \,M(x^{'})
+
\frac{1}{\check{p}}\,
\sum^{\check{p}}_\zs{j=l}\,\xi^2_\zs{j,p}
\,,
\end{equation}
where $M(x^{'})$ is defined in \eqref{sec:Pr.2}.
Furthermore, putting  
$
x_j^{''}= \check{p}^{-1/2}   \Chi_\zs{\{l \leqslant  j\leqslant \check{p}\}}$,
one can write the last term on the right hand side of \eqref{sec:Mo.1-0-1-04} as
$$ 
\frac{1}{\check{p}}\, \sum^{\check{p}}_\zs{j=l }\,\xi^2_\zs{j,p} = \frac{1}{\sqrt{\check{p}}} \,
B_\zs{2,Q}(x^{''})
 + \frac{1}{\check{p}} B_\zs{1,Q}(x^{'})  + \frac{(\check{p}- l+1) \sigma_\zs{Q }}{\check{p}}
\,,
$$
where  the functions $B_\zs{1,Q}$ and $B_\zs{2,Q}$ are given in \eqref{sec:Prsm.1}. 
Using Proposition~\ref{Le.sec:A.06-2}, Proposition~\ref{Le.sec:A.06-3} and Lemma \ref{Le.sec:A.2} , we come to the following upper bound 

$$
\E_\zs{Q}|\wh{\sigma}_\zs{n}-\sigma_\zs{Q}|
\le 
\frac{16 \Vert\dot{S}\Vert^{2} n  }{l p}
+\frac{2n}{p} \E_\zs{Q}\,| \mathbf{M}(x^{'})|
+\frac{\L_\zs{1,Q}}{p}
+\frac{\sqrt{\L_\zs{2,Q}}}{\sqrt{p}}
+\frac{\sigma_\zs{Q}(l -1)}{p}.
$$
In the same way as in
\eqref{sec:Pr.8}, we obtain
$$
\E_\zs{Q}\,|M(x^{'})| \le \left( \frac{\varkappa_\zs{Q}}{n}\,
\sum^{p}_\zs{j=l} \, \overline{\theta}^2_\zs{j,p}  \right) ^{1/2}
\le 
\frac{4(\varkappa_\zs{Q}\Vert\dot{S}\Vert^{2})^{1/2}}{l}.
$$

Taking into account that 
$\kappa_\zs{Q}\le (1+\check{\tau}\vert\rho\vert_\zs{*})\sigma_\zs{Q}$
and using the bounds
\eqref{sec:L_2-Upp_1-1} and
\eqref{sec:L_2-Upp_2}
we obtain the inequality
\eqref{sec:Si.3}. 
Hence Proposition \ref{Pr.sec:Si.1} holds true.
\fdem

 \bigskip
 \subsection{Proof of Theorem \ref{Th.sec:Mrs.1-10}}
This prof directly follows from Theorem \ref{Th.sec:Mrs.0} and Proposition \ref{Pr.sec:Si.1}.
\fdem

 \subsection{Proof of Theorem \ref{Th.sec:Ef.1}}
 
  First, we denote by $Q_\zs{0}$ the distribution of the noise
 \eqref{sec:Ex.1} and \eqref{sec:Mcs.1}
 with the parameter $\varrho_\zs{1}=\varsigma^{*}$, $\check{\varrho}=1$ and $\varrho_\zs{2}=0$, i.e., 
 the distribution for the ``signal + white noise''  
 model. So, we can estimate as below the robust risk
 $$
\cR^{*}_\zs{n}(\wt{S}_\zs{n},S)\ge
\cR_\zs{Q_\zs{0}}(\wt{S}_\zs{n},S)\,.
$$
Now, Theorem 6.1 from \cite{KonevPergamenshchikov2009b}
yields the lower bound \eqref{sec:Ef.4}.  Hence this finishes the proof. \fdem

\subsection{Proof of Proposition \ref{Th.sec:Ef.33}}

First, we note that in view of \eqref{sec:Mo.1}
one can represent the quadratic risk
for the empiric  norm $\Vert \cdot \Vert_\zs{p}$
as
$$
\E_\zs{Q}\,
\|\wh{S}_\zs{\lambda_\zs{0}}-S\|^{2}_\zs{p}
=
\frac{1}{n}\,
 \sum_{j=1}^{\check{p}}\,\lambda_\zs{0}^2(j)\,\E_\zs{Q}\,\xi^2_\zs{j,p}
+
\ov{\Theta}_\zs{p}\,,
$$
where $\ov{\Theta}_\zs{p}= \sum_{j=1}^{p}\,
\left(\theta_\zs{j,p}-\lambda_\zs{0}(j)\,\ov{\theta}_\zs{j,p}
\right)^2$. We put here $\lambda_\zs{0}(j)=0$ for $j>n$ if $p>n.$
The first term can be estimated  by the  bound
\eqref{sec:L_2-Upp_1-1}
 as
$$
\sup_\zs{Q\in\cQ_\zs{n}}\,
\E_\zs{Q}
\sum_{j=1}^{\check{p}}\,\lambda_\zs{0}^2(j)
\,\xi^{2}_\zs{j,p}\le 
\,\varsigma_\zs{*}
\,\sum_{j=1}^{n}\,\lambda_\zs{0}^2(j)
+\L_\zs{1,Q}\,.
$$
where $\L^{*}_\zs{1,n}=\sup_\zs{Q\in\cQ_\zs{n}}\,\L_\zs{1,Q}$.
Therefore, taking into account that $\upsilon_\zs{n}=n/\sigma^{*}$, we get
$$
\sup_\zs{Q\in\cQ_\zs{n}}\,
\E_\zs{Q}\,
\|\wh{S}_\zs{\lambda_\zs{0}}-S\|^{2}_\zs{p}
\,\le\,
\frac{1}{\upsilon_\zs{n}}
\sum_{j=1}^{n}\,\lambda_\zs{0}^2(j)
+\frac{\L^{*}_\zs{1,n}}{n}
+\ov{\Theta}_\zs{p}
\,.
$$
Note that
\begin{equation}\label{sec:Up.3}
\lim_{n\to\infty}\,
\frac{1}{\upsilon^{1/(2k+1)}_\zs{n}}\sum_{j=1}^{n}\,\lambda_\zs{0}^2(j)
=
\frac{2(\tau_\zs{k}\,r)^{1/(2k+1)}\,k^2}{(k+1)(2k+1)}
\,.
\end{equation}
\noindent Furthermore, by the Inequality \eqref{sec:Pr.5}
for any  $0<\wt{\ve}<1$ we get
\begin{equation}\label{sec:Up.3-0}
\ov{\Theta}_\zs{p}\le (1+\wt{\ve})\,\Theta_\zs{p}+
(1+\wt{\ve}^{-1})\,\sum^{p}_\zs{j=1}\,h^{2}_\zs{j,p}\,,
\end{equation}
where $\Theta_\zs{p}=
\sum^{p}_{j=1}\,(1-\lambda_\zs{0}(j))^2\,\theta^2_\zs{j,p}$.
In view of Definition \eqref{sec:Ga.2},
 we can represent this term as
$$
\Theta_\zs{p}=
\sum_{j=\iota_\zs{0}}^{[\omega_\zs{0}]}\,(1-\lambda_\zs{0}(j))^2\,\theta^2_\zs{j,p}
+
\sum_{j=[\omega_\zs{0}]+1}^{p}\,\theta^2_\zs{j,p}
:=\Theta_\zs{1,p}+\Theta_\zs{2,p}\,,
$$
where $\iota_\zs{0}=j_\zs{*}(\alpha_\zs{0})$,
$\omega_\zs{0}=\omega_\zs{\alpha_\zs{0}}=
\left(\tau_\zs{k} \l_\zs{0} \upsilon_\zs{n}\right)^{1/(2k+1)}$ and
$\l_\zs{0}=\left[r/ \varepsilon\right]\varepsilon$.
 Applying Lemma~\ref{Le.sec:A.4} yields
$$
\Theta_\zs{1,p}
\le (1+\wt{\varepsilon})\,
\sum_{j=l}^{[\omega_\zs{0}]}\,(1-\lambda_\zs{0}(j))^2\,\theta^{2}_\zs{j}
+4\pi^{2}r
(1+\wt{\varepsilon}^{-1})\,\omega^{3}_\zs{0}\,p^{-2}\,.
$$

\noindent 
Similarly, through Lemma~\ref{Le.sec:A.3} we have
$$
\Theta_\zs{2,p}\le (1+\wt{\varepsilon})
\sum_\zs{j\ge [\omega_\zs{0}]+1}\,\theta^{2}_\zs{j}
+(1+\wt{\varepsilon}^{-1})\,r\,p^{-2}\,.
$$

Hence,
$$
\Theta_\zs{p}\,
\le (1+\wt{\varepsilon})\,
\Theta^{*}_\zs{\iota_\zs{0}}
+
(1+\wt{\varepsilon}^{-1})\,
\left(4\pi^{2}r\omega^{3}_\zs{0}+r\right)\,p^{-2}\,,
$$
where $\Theta^{*}_\zs{l}
=\sum_\zs{j\ge l}\,(1-\lambda_\zs{0}(j))^2\,\theta^{2}_\zs{j}$.
Moreover,  note that
$$
\sup_\zs{S\in W^{1}_\zs{r}}\,
\max_\zs{1\le j\le p}\,
h^{2}_\zs{j,p}\le \,\|\dot{S}\|^{2}\,p^{-2}\,\le \,r\,p^{-2}\,.
$$
Moreover,  $W^{k}_\zs{r}\subseteq W^{2}_\zs{r}$  for any $k \ge 2$.
From here and Lemma~\ref{Le.sec:A.5}
 we get
$$
\sup_\zs{S\in W^{k}_\zs{r}}\,
\sum^{p}_\zs{j=1}\,h^{2}_\zs{j,p}\le r
\left(p^{-1}\,\Chi_\zs{\{k=1\}}+3p^{-2}\Chi_\zs{\{k\ge 2\}}\right)\,.
$$
Moreover, in view of 
 Condition $\H_\zs{5}$)
$$
\lim_\zs{n\to\infty}\,\upsilon^{2k/(2k+1)}_\zs{n}\,\left(p^{-1}\Chi_\zs{\{k=1\}}+\omega^{3}_\zs{0}p^{-2}\right)\,
=0\,.
$$
So,
$$
\limsup_\zs{n\to\infty}\,\upsilon^{2k/(2k+1)}_\zs{n}\,\sup_\zs{S\in W^{k}_\zs{r}}\,
\ov{\Theta}_\zs{p}\,
\le \,
\limsup_\zs{n\to\infty}\,\upsilon^{2k/(2k+1)}_\zs{n}\,\sup_\zs{S\in W^{k}_\zs{r}}\,
\Theta^{*}_\zs{\iota_\zs{0}}\,.
$$
To estimate the term $\Theta^{*}_\zs{\iota_\zs{0}}$
we set
$$
\U_\zs{n}= \upsilon^{2k/(2k+1)}_\zs{n} \sup_\zs{j\ge
\iota_\zs{0}}(1-\lambda_\zs{0}(j))^2/a_\zs{j}\,,
$$
where  the sequence $(a_\zs{j})_\zs{j\ge 1}$ is defined in \eqref{sec:Ef.2}.
This leads to the inequality
$$
\sup_\zs{S\in W^{1}_\zs{r}}\,
\upsilon^{2k/(2k+1)}_\zs{n}\,\Theta^{*}_\zs{\iota_\zs{0}}
\,
\le
\U_\zs{n}\,
\sum_\zs{j\ge 1}\,a_\zs{j}\,\theta^{2}_\zs{j}\,\le\,\U_\zs{n}\,r
\,.
$$
Taking into account that
$\lim_\zs{n\to\infty} t_\zs{0}=r$, we get
$$
\limsup_{n\to\infty}\,
\U_\zs{n}
\le\,
\pi^{-2k}\left(\tau_\zs{k}\,r \right)^{-2k/(2k+1)}
\,,
$$
where the coefficient $\tau_\zs{k}$ is given in \eqref{sec:Ga.2}.
This implies immediately that
\begin{equation}\label{sec:Up.4}
\limsup_\zs{n\to\infty}\,\upsilon^{2k/(2k+1)}_\zs{n}\,
\sup_\zs{S\in W^{k}_\zs{r}}\,
\overline{\Theta}_\zs{p}\,
\le
\frac{r^{1/(2k+1)}}{\pi^{2k}(\tau_\zs{k})^{2k/(2k+1)}}\,.
\end{equation}
Moreover, note that
$$
R^{*}_\zs{k}
=
\frac{2(\tau_\zs{k}\,r)^{1/(2k+1)}\,k^2}{(k+1)(2k+1)}
+
\frac{r^{1/(2k+1)}}{\pi^{2k}(\tau_\zs{k})^{2k/(2k+1)}}
\,.
$$
So, applying 
 \eqref{sec:Up.3} and \eqref{sec:Up.4}, yields
\begin{equation}\label{Up.5}
\lim_\zs{n\to\infty}
\upsilon^{2k/(2k+1)}_\zs{n}\,\sup_\zs{S\in W^{k}_\zs{r}}
\sup_\zs{Q\in\cQ_\zs{n}}\,
\E_\zs{Q}\,
\|\wh{S}_\zs{\lambda_\zs{0}}-S\|^{2}_\zs{p}
\le\,R^{*}_\zs{k}\,.
\end{equation}
Furthermore, Lemma~\ref{Le.sec:A.1-1} yields that for any 
$\wt{\varepsilon}>0$
$$
\sup_\zs{S\in W^{k}_\zs{r}}
\,\cR^{*}_\zs{n}(\wh{S}_\zs{\lambda_\zs{0}},S)
\le\,(1+\wt{\varepsilon})
\sup_\zs{S\in W^{k}_\zs{r}}
\sup_\zs{Q\in\cQ_\zs{n}}\,
\E_\zs{Q}\,
\|\wh{S}_\zs{\lambda_\zs{0}}-S\|^{2}_\zs{p}
+
\,(1+\wt{\varepsilon}^{-1})r\,p^{-2}\,.
$$
So, in view of Condition $\H_\zs{5}$),
we derive the desired inequality
$$
\lim_\zs{n\to\infty}\,
\upsilon^{2k/(2k+1)}_\zs{n}\,
\sup_\zs{S\in W^{k}_\zs{r}}
\,\cR^{*}_\zs{n}(\wh{S}_\zs{\lambda_\zs{0}},S)
\le\,R^{*}_\zs{k}\,.
$$
Hence we obtain Proposition \ref{Th.sec:Ef.33}.
 \fdem

\bigskip

{\bf Acknowledgments.}  
The last author is partially supported   by the Russian Federal Professor program  (Project No 1.472.2016/1.4, Ministry of Education and Science of the Russian Federation).

\bigskip

\renewcommand{\theequation}{A.\arabic{equation}}
\renewcommand{\thetheorem}{A.\arabic{theorem}}
\renewcommand{\thesubsection}{A.\arabic{subsection}}
\section{Appendix}\label{sec:A}
\setcounter{equation}{0}
\setcounter{theorem}{0}

\subsection{Property of the penalty term}

\begin{lemma}\label{Le.sec:A.1-06-11-01}
For any $n\ge\,1$ and $\lambda \in \Lambda$,
$$ 
P^{0}_n(\lambda) \leq 
 \cR_\zs{Q}(\wh{S}_\zs{\lambda},S)
+\frac{\L_\zs{1,Q}}{n}, 
$$
where the coefficient $P^{0}_n(\lambda)$ is defined in \eqref{sec:Pr.2} and the $\L_\zs{1,Q}$ is defined in\eqref{sec:L_2-Upp_1-1}.
\end{lemma}
\proof
 By the definition of $\Er(\lambda)$ in \eqref{sec:def-err.1}
  one has
$$ 
\Er(\lambda)\ge \sum_\zs{j=1}^{\check{p}} \left((\lambda(j)-1) \overline{\theta}_\zs{j,p}+ \frac{\lambda(j)}{n}\xi_\zs{j,p} \right)^2
\,. 
$$                                                                               
In view of Proposition~\ref{Le.sec:A.06-2} we obtain that
$$
\cR_\zs{Q}(\wh{S}_\zs{\lambda},S)=
\E_\zs{Q}\, \Er(\lambda) \ge\, \frac{1}{n}\sum_\zs{j=1}^{n} \lambda^2(j)  \E_\zs{Q}\,\xi^2_\zs{j,n} \ge\, P^{0}_n(\lambda)-\frac{\L_\zs{1,Q}}{n}
\,.
$$
Hence we otain Lemma \ref{Le.sec:A.1-06-11-01}.

\subsection{Properties of stochastic integrals \eqref{sec:In.2}}\label{sec:SI}

In this section we give some results of stochastic calculus for the process $(\xi_\zs{t})_\zs{t\ge\, 0}$ given in \eqref{sec:Ex.1}, needed all along this paper. As the process $\xi_\zs{t}$ is the combination of a L\'evy process and a semi-Markov process, these results are not standard and need to be provided.

\begin{lemma}\label{Le.sec:Smp.1} 
Assume that Conditions $\H_\zs{1})$--$\H_\zs{4})$ hold true. Then, for any $n\ge 1$ and for any non random function $f$ from $\L_\zs{2}[0,n],$ the stochastic integral \eqref{sec:In.2} exists and satisfies the properties \eqref{sec:In.3} with the coefficient $\varkappa_\zs{Q}$ given in 
\eqref{sec:In.3}.
\end{lemma}

\begin{lemma}\label{Le.sec:Smp.2} 
Let $f$ and $g$ be any non-random functions from 
$\L_\zs{2}[0,n]$ and $(I_\zs{t}(f))_\zs{t\ge\,0}$ be the process defined in \eqref{sec:In.2}. Then, for any $0 \leq t \leq n$,
\begin{equation}\label{sec:A.06-11-00-00}
\E\,I_\zs{t}(f) I_\zs{t}(g) = \bar{\varrho}
\,
(f,g)_\zs{t}\,
+ 
\varrho_\zs{3}^2
\,
(f,g\rho)_\zs{t}
\,, 
\end{equation}
where
$(f,g)_\zs{t}=\int_\zs{0}^{t} f(s)\,g(s) \d s$ and
 $\rho$ is the  density defined in \eqref{sec:Cns.1}.
\end{lemma}

\begin{lemma}\label{Le.sec:A.06-11-03} 
Let $f$ and $g$ be bounded functions defined on $[0,\infty) \times \bbr.$ Then, for any $k\ge 1,$
$$
\E \left( I_\zs{T_\zs{k-}} (f)\,I_\zs{T_\zs{k-}} (g) \mid \cG \right)= \bar{\varrho} 
(f\,,\,g)_\zs{T_\zs{k}}+ 
\varrho_\zs{3}^2 \sum_\zs{l=1}^{k-1}\, f(T_\zs{l})\,g(T_\zs{l}),
$$
where  $\cG$ is the $\sigma$-field generated by the sequence $(T_\zs{l})_\zs{l\ge 1}$, i.e., $\cG=\sigma\{T_\zs{l}\,,\,l\ge 1\}$.
\end{lemma}

\begin{lemma}\label{Le.sec:A.05-00} 
Assume that Conditions $\H_\zs{1})$--$\H_\zs{4})$ hold true. 
 Then, for any measurable bounded non-random functions $f$ and $g$, one has
$$ 
\E \int_\zs{0}^{n} I^2_\zs{t-}(f)  I_\zs{t-}(g) g(t) d\xi_\zs{t} = 0.
$$
\end{lemma}

Lemmas  \ref{Le.sec:Smp.1} -- \ref{Le.sec:A.05-00}
are proved in \cite{VladSlimSerge2016}.


\subsection{Properties of the Fourier coefficients}\label{sec:Coef}

\begin{lemma}\label{Le.sec:A.1-1}
Let $f$ be an absolutely continuous function, $f: [0,1]\to\bbr,$ with
$\|\dot{f}\|<\infty$ and $g$ be a simple function, $g: [0,1]\to\bbr$ of the form
$
g(t)=\sum_{j=1}^p\,c_\zs{j}\,\chi_{(t_{j-1},t_j]}(t),$
where $c_\zs{j}$ are some constants. Then for any $\varepsilon>0,$  the function $\Delta=f-g$
satisfies the following inequalities
$$
\|\Delta\|^{2}\le (1+\wt{\varepsilon})\|\Delta\|^{2}_\zs{p}
+
(1+\wt{\varepsilon}^{-1})\frac{\|\dot{f}\|^{2}}{p^{2}}\,,
\quad
\|\Delta\|^{2}_\zs{p}\le (1+\wt{\varepsilon})\|\Delta\|^{2}
+
(1+\wt{\varepsilon}^{-1})\frac{\|\dot{f}\|^{2}}{p^{2}}
\,.
$$
\end{lemma}

\begin{lemma}\label{Le.sec:A.2} 
Let the function $S(t)$ in \eqref{sec:In.1} be absolutly continuous and have an absolutely integrable derivative. Then the coefficients $(\overline{\theta}_\zs{j,p})_{1\leqslant j \leqslant p}$ defined in \eqref{sec:In.08} satisfy the inequalities
\begin{equation}
\label{sec:A.11-00-1}
\vert \overline{\theta}_\zs{1,p} \vert \leqslant \Vert S\Vert_1  \quad\mbox{and}\quad\
 \max_\zs{2\leqslant j \leqslant p} j \vert \overline{\theta}_\zs{j,p} \vert
\leqslant 2 \sqrt{2} \Vert \dot{S} \Vert_\zs{1}\,.
\end{equation}
\end{lemma}

\begin{lemma}\label{Le.sec:A.3}
For any $p\ge 2$, $1\le N\le p$ and $r>0$,  
the coefficients $(\theta_\zs{j,p})_\zs{1\le j\le p}$
of functions $S$ from
the class $W^{1}_\zs{r}$ satisfy, for any $\wt{\varepsilon}>0$,
the following inequality
\begin{equation}\label{sec:A.2}
\sum^{p}_\zs{j=N}
\theta^{2}_\zs{j,p}
\,
\le\,(1+\wt{\varepsilon})
\,\sum_\zs{j\ge N}\,\theta^{2}_\zs{j}
\,
+(1+\wt{\varepsilon}^{-1})\,r\,p^{-2}
\,.
\end{equation}
\end{lemma}

\begin{lemma}\label{Le.sec:A.4}
For any $p\ge 2$ and $r>0$,   the coefficients $(\theta_\zs{j,p})_\zs{1\le j\le p}$
of functions $S$ from
the class $W^{1}_\zs{r}$ satisfy
the following inequality
\begin{equation}\label{sec:A.3}
\max_\zs{1\le j\le p}
\,\sup_\zs{S\in W^{1}_\zs{r}}
\left(
|\theta_\zs{j,p}
-
\theta_\zs{j}|
-2\pi \sqrt{r} \,j\,p^{-1}
\right)
\,
\le 0\,.
\end{equation}
\end{lemma}

\begin{lemma}\label{Le.sec:A.5}
For any $p\ge 2$ and $r>0$
the correction coefficients 
$\left(h_\zs{j,p}\right)_\zs{1\le j\le p}$
 for the functions $S$ from
the class $W^{2}_\zs{r}$ satisfy
the following inequality
\begin{equation}\label{sec:A.4}
\sup_\zs{S\in W^{2}_\zs{r}}
\sum^{p}_\zs{j=1}
h^{2}_\zs{j,p}
\le \,3r\,p^{-2}
\,.
\end{equation}
\end{lemma}

\noindent 
Lemmas 
\ref{Le.sec:A.1-1}
--
\ref{Le.sec:A.5}
are proven in 
\cite{KonevPergamenshchikov2015}.


\medskip

\medskip

\newpage

\end{document}